\documentclass[twoside]{article}
\usepackage{amsmath,amsthm,amssymb,amscd}
\usepackage[all]{xy}
\theoremstyle{plain}
 
\newtheorem{lem}{Lemma}[section]
\newtheorem{cor}{Corollary}[section]

\newtheorem{thm}{Theorem}[section]

\newcommand{\bb}[1]{\mbox{$\mathbb{#1}$}}
\theoremstyle{definition}
\newtheorem{ex}{Example}[section]
\theoremstyle{remark}
\newtheorem{rem}{Remark}[section]

\title{Hirzebruch classes and motivic Chern classes\\ for singular 
spaces\footnote{This is a completely new and improved version of the paper math.AG/0405412.}}
\author{Jean-Paul Brasselet \and J\"{o}rg Sch\"{u}rmann \and Shoji Yokura}
\date{{\em Dedicated to the memory of Shiing-shen Chern\\ and to Friedrich Hirzebruch}}

\begin{document}

\maketitle
\bibliographystyle{plain}

\begin{abstract}
In this paper we study some new theories of characteristic homology classes
of singular complex algebraic (or compactifiable analytic) spaces.

We introduce a {\em motivic Chern class transformation\/} $mC_{*}: K_{0}(var/X)$ $\to
G_{0}(X)\otimes \bb{Z}[y]$, which generalizes the total $\lambda$-class
$\lambda_{y}(T^{*}X)$ of the cotangent bundle to singular spaces. Here 
$K_{0}(var/X)$ is the relative Grothendieck group of complex algebraic varieties over $X$ as introduced and 
studied by Looijenga and Bittner in relation to motivic integration, and $G_{0}(X)$ is the Grothendieck group 
of coherent sheaves of ${\cal O}_{X}$-modules. 
A first construction of
$mC_{*}$ is based on resolution of singularities and a suitable
``blow-up'' relation. In the (complex) algebraic context this ``blow-up'' relation follows  from work of Du Bois, based on Deligne's mixed Hodge theory.
Other approaches by work of Guill\'en and Navarro Aznar (using ``only'' resolution of singularities) or Looijenga and Bittner (using the ``weak factorization theorem'')
also apply to the compactifiable complex analytic context.
A second more functorial construction of $mC_{*}$ is based on some results from the theory of algebraic mixed Hodge modules due to M.Saito.  

We define a natural
transformation $T_{y*}: K_{0}(var/X) \to H_{*}(X)\otimes \bb{Q}[y]$ commuting
with proper pushdown, which generalizes the corresponding {\em Hirzebruch
characteristic\/}. $T_{y*}$ is a homology class
version of the motivic measure corresponding to a suitable specialization
of the well known Hodge polynomial. This transformation unifies
the {\em Chern class transformation\/} of MacPherson and Schwartz (for $y=-1$), the 
{\em Todd class transformation\/} in the singular Riemann-Roch theorem of
Baum-Fulton-MacPherson (for $y=0$) and the {\em $L$-class transformation\/} of Cappell-Shaneson (for $y=1$).  

In the simplest case of a normal Gorenstein variety with ``canonical singularities'' we also explain a relation among the
``stringy version'' of our characteristic classes, the {\em elliptic class\/} of
Borisov-Libgober and the {\em stringy Chern classes\/} of Aluffi and 
De Fernex-Lupercio-Nevins-Uribe.

Moreover, all our results can be extended to varieties over a base field $k$ of characteristic $0$.
\end{abstract}

\section*{Introduction}

In this paper we study some new theories of characteristic homology classes
of a singular algebraic variety $X$ defined over a base field $k$ of characteristic $0$.

Let $K_{0}(var/X)$ be the {\em relative Grothendieck group of
algebraic varieties over $X$\/} as introduced and studied by Looijenga \cite{Lo}
and Bittner \cite{Bi} in relation to motivic integration.
Here a variety is a seperated scheme of finite type over $spec(k)$.
$K_{0}(var/X)$ is the quotient of the free abelian group of isomorphism
classes of algebraic morphisms $Y\to X$ to $X$,  modulo the ``additivity''
relation generated by  
\begin{equation} \label{eq:add}
[Y\to X] = [Z\to Y \to X] + [Y\backslash Z \to Y \to X] 
\tag{add}
\end{equation}
for $Z\subset Y$ a closed algebraic subvariety of $Y$.
Taking $Z=Y_{red}$ we see that these classes depend only on the underlying
reduced spaces.
By resolution of singularities, $K_{0}(var/X)$ is generated by
classes $[Y\to X]$ with $Y$ smooth, pure dimensional and proper
over $X$. Moreover, for any morphism $f: X'\to X$ we get a functorial pushdown:
\[f_{!}: K_{0}(var/X')\to K_{0}(var/X);\: [h:Z\to X']\mapsto [f\circ h: Z\to X] \:.\]

We introduce characteristic homology class transformations $mC_{*}$ and $T_{y*}$
on $K_{0}(var/X)$ related by the following commutative
diagram:
\begin{equation} \label{eq:mcTy} \begin{CD} 
K_{0}(var/X) @= K_{0}(var/X)  \\ 
@VV mC_{*} V  @VV  T_{y*} V\\
G_{0}(X)\otimes \bb{Z}[y] 
@> td_{(1+y)} >> H_{*}(X)\otimes \bb{Q}[y,(1+y)^{-1}]
\supset  H_{*}(X)\otimes \bb{Q}[y]  \:,
\end{CD} \end{equation}
where $H_{*}(X)$ denotes the Chow group $A_{*}(X)$, or the Borel-Moore
homology $H^{BM}_{2*}(X,\bb{Z})$ (in even degrees) in the complex context.
Moreover, $G_{0}(X)$ is the Grothendieck group 
of coherent sheaves of ${\cal O}_{X}$-modules. 
Here $td_{(1+y)}$ is a generalization given by Yokura \cite{Y4}
of the {\em Todd class transformation\/} $td_{*}$ used in the singular Riemann-Roch
theorem of Baum-Fulton-MacPherson \cite{BFM,FM} (for Borel-Moore homology)
or Fulton \cite{Ful} (for Chow groups).\\

Our transformations $mC_{*}$ and $T_{y*}$ commute with proper pushdown 
and are uniquely determined by this property
together with the following normalization conditions for $X$ smooth and
pure $d$-dimensional:
\begin{itemize}
\item $mC_{*}([id_{X}])= \sum_{i=0}^{d} \; [\Lambda^{i} T^{*}X]\cdot y^{i}
= \lambda_{y}([T^{*}X])\cap [{\cal O}_{X}]$,
with $\lambda_{y}$ the total $\lambda$-class (compare \cite{FL}).
\item $T_{y*}([id_{X}])= T^{*}_{y}(TX) \cap [X]$,
with $T^{*}_{y}$ the corresponding {\em Hirzebruch characteristic class\/} introduced in \cite{Hi} (compare with Section \ref{sec:gHHR}).
\end{itemize}

Moreover, the transformation $T_{y*}$ fits into the commutative diagram
\begin{equation} \label{eq:real1} \begin{CD}
F(X) @<e<<  K_{0}(var/X) @> mC_{0} >> G_{0}(X) \\
@V c_*\otimes\bb{Q} VV   @V T_{y*} VV   @VV td_{*} V \\
H_{*}(X)\otimes\bb{Q}  @<y=-1<< H_{*}(X)\otimes\bb{Q}[y] @> y=0 >>
H_{*}(X)\otimes\bb{Q} \:,
\end{CD} \end{equation}
with $c_{*}\otimes\bb{Q}$ the (rationalized) {\em Chern-Schwartz-MacPherson transformation\/} 
on the group $F(X)$ of algebraically constructible functions \cite{Schwa,M1,BrS,Ken}.
Here $mC_{0}$ is the degree zero component of our {\em motivic
Chern class transformation\/} $mC_{*}$, and $e$ is simply given by
\begin{equation} \label{eq:e}e([f: Y\to X]):=f_{!}1_{Y} \in F(X)\:,
\end{equation}
i.e., by taking fiberwise the (topological) Euler characteristic with compact support.
Note that the homomorphisms $e$ and $mC_{0}$ are surjective.

\begin{rem} \label{rem:c*alg}
Note that in the algebraic context the {\em Chern-Schwartz-MacPherson transformation\/} $c_{*}$
of \cite{Ken} is defined only on spaces $X$ embeddable into a smooth space
(e.g. for quasi-projective varieties). But using the technique of ``Chow envelopes''
as in \cite[sec.18.3]{Ful}, this transformation can uniquely be extended to 
all (reduced) seperated schemes of finite type over $spec(k)$.
\end{rem}

For $X$ a compact complex algebraic variety, we can also construct the 
following commutative diagram of natural transformations:
\begin{equation} \label{eq:real2} \begin{CD}
K_{0}(var/X) @> sd >>  \Omega(X)\\
@V T_{y*} VV   @VV L_{*} V \\
H_{2*}(X,\bb{Q})[y] @> y=1 >>
H_{2*}(X,\bb{Q}) \:,
\end{CD} \end{equation}
with $L_{*}$ the {\em homology $L$-class transformation\/} of Cappell-Shaneson \cite{CS1}
(as reformulated by Yokura \cite{Y1}). Here $\Omega(X)$ is the abelian group
of cobordism classes of selfdual constructible complexes. To make $\Omega(\cdot)$
covariant functorial (correcting \cite{Y1}), we use the definition
of a ``cobordism'' given by Youssin \cite{You} in the more general context of
triangulated categories with duality.
Note that the {\em homology $L$-class transformation\/} of Cappell-Shaneson \cite{CS1}
is defined only for compact spaces and takes values in usual homology,
since its definition is based on a corresponding signature invariant together
with the Thom-Pontrjagin construction.
Also $sd$ does not need to be surjective.\\

So $T_{y*}$ unifies the (rationalized) {\em Chern-Schwartz-MacPherson
transformation\/} $c_{*}\otimes\bb{Q}$, the {\em Todd transformation\/} $td_{*}$ of
Baum-Fulton-MacPherson (for Borel-Moore homology) or Fulton (for Chow groups)
and the {\em $L$-class transformation\/} of Cappell-Shaneson.
However, we point out that our approach does not give a new proof of the existence
of these (now) classical transformations, since our approach only gives the
corresponding transformations
 on the relative Grothendieck group $K_{0}(var/X)$.
Moreover, the {\em Hirzebruch class\/} $T_{y*}(X):=T_{y*}([id_{X}])$ of the singular space
$X$ specializes, for $y=-1$, to the (rationalized) {\em Chern-Schwartz-MacPherson class\/}
$c_{*}(X)\otimes\bb{Q}:=c_{*}(1_{X})\otimes\bb{Q}$ of $X$, since $e([id_{X}])=1_{X}$.
But, in general we have: for $y=0$
\[mC_{0}([id_{X}])\neq [{\cal O}_{X}] \quad \text{and} \quad T_{0*}(X) \neq 
td_{*}(X):=td_{*}([{\cal O}_{X}]) \:,\]
and similarly for $y=1$
\[sd([id_{X}])\neq [{\cal IC}_{X}] \quad\text{and} \quad T_{1*}(X) \neq 
L_{*}(X):=L_{*}([{\cal IC}_{X}]) \:,\]
with ${\cal IC}_{X}$ the middle {\em intersection cohomology complex\/} of Goresky-MacPherson
\cite{GM}. This means that our approach based on the ``additivity'' property picks out new
distinguished elements 
\[mC_{0}([id_{X}])\in G_{0}(X) \quad \text{and} \quad
sd([id_{X}]) \in \Omega(X)\:.\]
On the positive side, we can show $mC_{0}([id_{X}]) = [{\cal O}_{X}]$ if $X$ has 
at most ``Du Bois singularities'' (e.g. rational singularities \cite{Kov, Sai5}). 
And similarly, we {\em conjecture\/} $sd([id_{X}])= [{\cal IC}_{X}]$ for $X$ a {\em rational homology manifold\/}.\\
 
\begin{rem} \label{rem:not}
Maybe here is the right place to explain our notion {\em motivic Chern class
transformation\/} for $mC_{*}$. Of course the notion
``motivic (dual) $\lambda$-class transformation'' would also be possible
by the corresponding normalization condition for $X$ smooth.
But we understand our transformations $T_{y*}$ and $mC_{*}$ by (\ref{eq:mcTy}) and
(\ref{eq:real1})
as natural ``motivic liftings'' of the {\em Chern-Schwartz-MacPherson transformation\/}
$c_{*}$, with $T_{-1*}([id_{X}])=c_{*}(X)\otimes\bb{Q}$ for any singular $X$.
\end{rem}

Our first construction of the motivic Chern class transformation $mC_{*}$ is based
on the following simple description of $K_{0}(var/X)$ in terms of proper morphisms
$[X'\to X]$. Let $Iso^{pr}(var/X)$ be the free abelian group on isomorphism classes
of proper morphisms $[X'\to X]$. Then one gets a canonical quotient map
\[Iso^{pr}(var/X)\to K_{0}(var/X)\]
(commuting with proper pushdown), which is an epimorphism of groups by 
``additivity''.

\begin{lem} \label{lem:iso1}
$K_{0}(var/X)$ is isomorphic to the quotient of $Iso^{pr}(var/X)$ modulo the
``acyclicity'' relation
\begin{equation} \label{eq:ac}
[\emptyset \to X]=0 \quad \text{and} \quad
[\tilde{X}'\to X] - [\tilde{Z}'\to X]= [X'\to X] - [Z'\to X] \:,\tag{ac}
\end{equation}
for any cartesian diagram
\begin{displaymath} \begin{CD}
\tilde{Z}' @>>> \tilde{X}' \\
@VVV  @VV q V \\
Z' @> i >> X' @>>> X \:,
\end{CD} \end{displaymath}
with $q$ proper, $i$ a closed embedding and $q: \tilde{X}'\backslash \tilde{Z}'
\to X' \backslash Z'$ an isomorphism.
\end{lem}

For the base field $k=\bb{C}$ and $q: X'\to X$ proper we now consider the {\em filtered
Du Bois complex\/} $(\underline{\Omega}^{*}_{X'},F)$ of $X'$ as introduced in \cite{DB}. 
This is a filtered complex of ${\cal O}_{X'}$-modules with maps differential operators
of order at most one, whose graded pieces 
\[gr_{F}^{p}(\underline{\Omega}^{*}_{X'}) \in D^{b}_{coh}(X')\]
are bounded complexes of ${\cal O}_{X'}$-modules with coherent cohomology.
Here the filtration is decreasing.
This filtered Du Bois complex is unique up to isomorphism in a suitable 
derived category $D_{diff}(X')$.
Here the assumption $k=\bb{C}$ is used, since the proof of \cite{DB} depends
on Deligne's theory of mixed Hodge structures.
In particular,
\[[gr_{F}^{p}(\underline{\Omega}^{*}_{X'})]:= \sum \nolimits_{i}\:(-1)^{i} 
H^{i}(gr_{F}^{p}(\underline{\Omega}^{*}_{X'})) \in G_{0}(X')\]
is well defined. Moreover, $gr_{F}^{p}$ commutes with proper pushdown 
\cite[prop.1.3]{DB}. By the construction of \cite{DB} one has 
\[gr_{F}^{p}(\underline{\Omega}^{*}_{X'})\simeq 0 \quad \text{for $p<0$ or $p>dim(X')$}\:,\]
together with a filtration-preserving map from the algebraic De Rham complex
(with the ``stupied'' filtration, compare Section \ref{sec:Hodge})
\[can: (DR({\cal O}_{X'}),F) \to (\underline{\Omega}^{*}_{X'},F)\:,\]
which is a filtered quasi-isomorphism for $X'$ smooth. Then the transformation
\begin{equation} \label{eq:mC1}
\begin{split}
mC_{*}: Iso^{pr}(var/X)\to G_{0}(X)\otimes \bb{Z}[y]\:; \\
[q: X'\to X] \mapsto \sum_{p \geq 0} \: q_{*}[gr_{F}^{p}(\underline{\Omega}^{*}_{X'})]\cdot (-y)^{p}
\end{split} \end{equation}
satisfies, by \cite[prop.1.3,prop.4.11]{DB}, the ``acyclicity'' relation (\ref{eq:ac})
and therefore, by lemma \ref{lem:iso1}, induces our motivic Chern class
transformation $mC_{*}$. In particular 
\[mC_{p}([id_{X}])= [gr_{F}^{p}(\underline{\Omega}^{*}_{X})[-p]]=
(-1)^{p}\cdot  [gr_{F}^{p}(\underline{\Omega}^{*}_{X})] \in G_{0}(X)\]
and by definition $X$ has at most ``Du Bois singularities'' if
\[can: {\cal O}_{X}= gr_{F}^{0}(DR({\cal O}_{X})) \to
gr_{F}^{0}(\underline{\Omega}^{*}_{X})\]
is a quasi-isomorphism.  \\

The more recent results of Guill\'en and Navarro Aznar \cite{GNA}
allow a construction of a Du Bois complex $\underline{\Omega}^{*}_{X}$
as before without the use of (mixed) Hodge theory. Instead, ``only'' 
Hironaka's resolution of singularities is used \cite[(2.1.2)]{GNA},
together with a corresponding ``Chow lemma'' \cite[(2.1.3)]{GNA}
for the comparison of different resolutions. So by this approach
we can define $mC_{*}$ as before in the algebraic context over any base field
$k$ of characteristic $0$.\\

Similarly it applies to the context of
compactifiable complex analytic spaces \cite[thm.4.1, sec.4.5-4.7]{GNA},
using Hironaka's resolution of singularities for analytic spaces
together with an ``analytic Chow lemma'' \cite[lem. on p.69]{GNA}.
Here we work in the category $AN_{\infty}$ of {\em compactifiable complex analytic spaces\/}
$X= \bar{X}\backslash \partial X$, with $\partial X$ a closed analytic subspace of
the compact complex analytic space $\bar{X}$.
Then we fix an equivalence class $\{(\bar{X},\partial X)\}$
of bimeromorphic equivalent compactifications, and consider only similar equivalence
classes of complex analytic morphism $f: X' \to X$ with a holomorphic extension
$\bar{f}: \bar{X}'\to \bar{X}$ (compare \cite[sec.4.5-4.7]{GNA}).
And for a compact complex analytic space $X$ we introduce
the {\em analytic relative Grothendieck group\/} $K_{0}(an/X)$ as the quotient of the free abelian group of isomorphism
classes of compactifiable analytic morphisms $Y\to X$ to $X$, modulo the ``additivity''
relation (\ref{eq:add})
for $Z\subset Y$ a compactifiable inclusion of a closed analytic subspace of $Y$.
Then lemma \ref{lem:iso1} holds also in this analytic context, with $Iso^{pr}(an/X)$
the free abelian group on isomorphism classes
of proper analytic morphisms $[X'\to X]$ (i.e. with $X'$ compact).\\

Note that the Chern-Schwartz-MacPherson class transformation $c_{*}$ and the 
$L$-class transformation of Cappell-Shaneson are also defined for compact spaces in this complex analytic context. 
The Todd transformation
$td_{*}: G_{0}(X)\to H_{2*}(X,\bb{Q})$ can be defined as
the composition
\begin{equation} \label{eq:toddan} \begin{CD}
G_{0}(X) @> \alpha >> K^{top}_{0}(X) @> td_{*} >> H_{2*}(X,\bb{Q}) \:,
\end{CD} \end{equation}
with $\alpha$ the {\em Riemann-Roch transformation\/} to (periodic) topological K-homology
(in even degrees) constructed by Levy \cite{Levy} (generalizing the
corresponding transformation of Baum-Fulton-MacPherson \cite{BFM2} for the
quasi-projective complex algebraic context) and the {\em topological Todd transformation\/} $td_{*}$ constructed in \cite{BFM2,FM}. If we restrict ourselve to a compact complex {\em manifold\/}
$X$, then we can also apply the {\em Riemann-Roch transformation to Hodge cohomology\/}
\begin{equation} \label{eq:toddhodge} \begin{CD}
G_{0}(X) @> \tau >>  \bigoplus_{k \geq 0} \: H^{k}(X,\Omega_{X}^{k}) 
\end{CD} \end{equation}
constructed by O'Brian-Toledo-Tong \cite{OTT}. In any case we also get the
commutative diagrams (\ref{eq:mcTy}), (\ref{eq:real1}) and (\ref{eq:real2})
for a compact complex analytic space $X$.\\

The simplest approach to our characteristic classes and the corresponding
natural transformations comes from a description of $K_{0}(var/X)$ 
(or $K_{0}(an/X)$) in terms of a ``blow-up'' relation for {\em smooth spaces\/}
mapping {\em properly\/} to $X$, which is due to Looijenga \cite{Lo} and Bittner \cite{Bi}.
This presentation is an easy application of the deep ``weak factorization theorem'' of \cite{AKMW,W}.\\

Let $Iso^{pr}(sm/X)$ be the free abelian group on isomorphism classes
of proper morphisms $[X'\to X]$, with $X'$ smooth
(and pure dimensional and/or quasi-projective if one wants). Then one gets a canonical quotient map
\[Iso^{pr}(sm/X)\to K_{0}(var/X)\]
(commuting with proper pushdown), which is an epimorphism of groups by ``additivity''
and Hironaka's resolution of singularities. Then we have by \cite[thm.5.1]{Bi}
the following basic result:

\begin{thm}[Bittner] \label{thm:blowup}
The group $K_{0}(var/X)$ is isomorphic to the quotient of $Iso^{pr}(sm/X)$ modulo the
``blow-up'' relation
\begin{equation} \label{eq:bl}
[\emptyset \to X]=0 \quad \text{and} \quad
[Bl_{Y}X'\to X] - [E\to X]= [X'\to X] - [Y\to X] \:,\tag{bl}
\end{equation}
for any cartesian diagram
\begin{displaymath} \begin{CD}
E @> i'>> Bl_{Y}X' \\
@VV q' V  @VV q V \\
Y @> i >> X' @> f >> X \:,
\end{CD} \end{displaymath}
with $i$ a closed embedding of smooth (pure dimensional) spaces and 
$f:X'\to X$ proper. Here $Bl_{Y}X'\to X'$ is the blow-up of $X'$ along 
$Y$ with exceptional divisor $E$. Note that all these spaces over $X$ are
also smooth (and pure dimensional and/or quasi-projective).
\hfill $\Box$
\end{thm}

\begin{rem} \label{rem:blowupan}
Since ``resolution of singularities''  and the ``weak factorization theorem''
can also be used in the complex analytic context, the simple proof of
\cite[thm.5.1]{Bi} implies the counterpart of Theorem \ref{thm:blowup}
also for $K_{0}(an/X)$ with $X$ a compact complex analytic space.
Of course here we cannot require $X'$ to be quasi-projective!
\end{rem}

A big advantage of Theorem \ref{thm:blowup} is the fact
that the ``blow-up'' relation (\ref{eq:bl}) is a statement for
{\em smooth\/} (pure dimensional) spaces mapping {\em properly\/} to $X$. Therefore all
our transformations are determined by ``functoriality'' and ``normalization''.
Then we only have to check the ``blow-up'' relation (\ref{eq:bl}) in the
following form (and similarly for the complex analytic context):

\begin{cor} \label{cor:bl}
Let $B_{*}: var/k \to group$ be a functor from the category $var/k$ of
(reduced) seperated schemes of finite type over $spec(k)$ to the category of abelian groups,
which is covariantly functorial
for proper morphism, with $B_{*}(\emptyset):=\{0\}$. Assume we can associate to any 
(quasi-projective) smooth space $X\in ob(var/k)$ 
(of pure dimension) a distinguished element $d_{X}\in B_{*}(X)$ such that
\begin{enumerate}
\item $h_{*}(d_{X'})=d_{X}$ for any isomorphism $h: X'\to X$.  
\item $q_{*}(d_{Bl_{Y}X})-i_{*}q'_{*}(d_{E}) = d_{X}-i_{*}(d_{Y}) \in B_{*}(X)$
for any cartesian blow-up diagram as in theorem \ref{thm:blowup} with $f=id_{X}$.
\end{enumerate}
Then there is by (1) a unique group homomorphism $\Phi : Iso^{pro}(sm/X) \to B_{*}(X)$
satisfying the ``normalization'' condition $\Phi([f:X'\to X])=f_{*}(d_{X'})$.
By (2) and {\em functoriality\/} this satisfies the ``blow-up'' relation of Theorem \ref{thm:blowup}
so that there is a unique induced group homomorphism $\Phi : K_{0}(var/X) \to B_{*}(X)$
commuting with proper pushdown and satisfying the ``normalization'' condition
$\Phi([id_{X}])=d_{X}$ for $X$ (quasi-projective) smooth (and pure dimensional).
\hfill $\Box$
\end{cor}

As a first example, we get the existence of a unique transformation
\[mC_{0}: K_{0}(var/X) \to G_{0}(X)\] 
with $mC_{0}([id_{X}])=[{\cal O}_{X}]$
for $X$ smooth (and pure dimensional). Note that the ``blow-up'' relation (2)
of Corollary \ref{cor:bl} follows from the well known relations
(compare \cite[(R5) on p.106, prop.4.1 on p.169]{FL}):
\[q'_{*}[{\cal O}_{E}]= [{\cal O}_{Y}] \quad \text{and} \quad
q_{*}[{\cal O}_{Bl_{Y}X}]= [{\cal O}_{X}] \:.\]

Similarly, the existence of a unique transformation
\[mC_{*}: K_{0}(var/X) \to G_{0}(X)\otimes \bb{Z}[y]\] 
with $mC_{*}([id_{X}])=
\sum_{i=0}^{d} \; [\Lambda^{i} T^{*}X]\cdot y^{i}
= \lambda_{y}([T^{*}X])\cap [{\cal O}_{X}]$
for $X$ smooth and pure $d$-dimensional follows directly from
\cite[IV.1.2.1]{Gr} or \cite[prop.3.3]{GNA} (where the last reference
also applies to the complex analytic context).\\

From the ``blow-up'' relation for $\lambda_{y}([T^{*}X])\cap [{\cal O}_{X}]$
we can deduce the ``blow-up'' relation for the {\em Hirzebruch 
class\/} $T^{*}_{y}(TX)\cap [X] \in H_{*}(X)\otimes \bb{Q}[y]$ in the context of 
quasi-projective smooth spaces
by using ``only'' the classical {\em Grothendieck-Riemann-Roch theorem\/} for 
quasi-projective smooth spaces as in \cite[thm.15.2]{Ful}
(instead of the (modified) singular Riemann-Roch transformation as in
(\ref{eq:mcTy})).\\

Finally, in Section \ref{sec:Lclass} we deduce  from Corollary \ref{cor:bl} the existence
of the unique {\em selfduality transformation\/} $sd: K_{0}(var/X) \to \Omega(X)$ with
$sd([id_{X}])= [\bb{Q}_{X}[n]]$ for $X$ smooth and pure $n$-dimensional.\\

In the final section we explain, for $k\subset\bb{C}$, the most powerful and functorial construction 
of the motivic Chern class transformation $mC_{*}$.
This was in the beginning our original approach. It is based on some deep results from the theory of {\em algebraic mixed Hodge modules\/} 
due to M.Saito \cite{Sai1}-\cite{Sai6}, which imply the existence of the two natural transformations 
\begin{displaymath} \begin{CD}
K_{0}(var/X) @> mH >> K_{0}(MHM(X/k))
@> gr^{F}_{-*}DR >>
H_{*}(X)\otimes \bb{Z}[y,y^{-1}] \,
\end{CD} \end{displaymath}
whose composition is our {\em motivic transformation\/} $mC_{*}$. 
Here $K_{0}(MHM(X/k))$ is the Grothendieck group of the abelian category
$MHM(X/k)$ of mixed Hodge modules on $X$, and $mH$ is defined by 
\[mH([f: Y\to X]):= [f_{!}\bb{Q}^{H}_{Y}]\in  K_{0}(MHM(X/k)) \:.\] 
Here $\bb{Q}^{H}_{Y}$ is in
some sense the ``constant Hodge module''  on $Y$. Similarly, $gr^{F}_{-*}DR$
comes form a suitable {\em filtered de Rham complex\/} of the filtered holonomic
D-module underlying a mixed Hodge module.\\

In this paper we focus on the construction of our homology class
transformations $mC_{*}$ and $T_{y*}$, together with the {\em unification\/}
of the transformations $c_*\otimes\bb{Q},\: td_{*}$ and $L_{*}$.

In a sequel to this paper we will show that $mC_{*}, T_{y*}$
and the Chern-Schwartz-MacPherson transformation $c_{*}$ are in fact {\em universal additive characteristic classes\/} that can be defined as natural transformations 
on the relative Grothendieck group $K_{0}(var/X)$.
The appropiate notion of a general (co)homology theory with Chern class
operators  was already introduced by Levine and Morel \cite{LM} in the algebraic context
as an {\em oriented cohomology theory\/} or an {\em oriented Borel-Moore weak homology theory\/}. And two of their main results describe 
\begin{itemize}
\item algebraic $K$-theory $K^{0}(X)=G_{0}(X)$ as the universal {\em oriented cohomology theory\/}
with a {\em multiplicative formal group law\/} on the category of smooth quasi-projective spaces, and
\item Chow groups $A_{*}(X)$ as the universal {\em oriented Borel-Moore weak homology theory\/}
with an {\em additive formal group law\/} on $var/k$.
\end{itemize}

Here we only sketch our main idea in the context of a classical {\em characteristic
cohomology class transformation\/} $cl^{*}: K^{0}(X)\to H^{2*}(X,R)$ for complex vector
bundles (with $R$ a $\bb{Q}$-algebra), which is {\em multiplicative\/}
\[cl^{*}(V \times W) = cl^{*}(V) \times cl^{*}(W) \]
and {\em normalized\/}: $cl^{0}(V)=1 \in H^{0}(\cdot,R)$.
And we assume that the corresponding transformation (for the base field $k=\bb{C}$)
\[\Phi_{cl} : Iso^{pro}(sm/X) \to H^{BM}_{2*}(X,R) \: ;
[f:X'\to X] \mapsto f_{*}(cl^{*}(TX')\cap [X']) \]
satisfies the ``blow-up'' relation (2) of Corollary \ref{cor:bl}
in the special case $X=\{pt\}$ a point. Then $\Phi=\Phi_{cl}$ commutes with exterior products
(by the multiplicativity of $cl^{*}$) and factorizes as ring homomorphisms

\begin{equation} \label{eq:genus} \begin{CD} 
Iso^{pro}(sm/\{pt\}) @>>> \Omega^{U}_{*}\otimes\bb{Q}  \\ 
@VVV  @VV \Phi V\\
K_{0}(var/\{pt\}) 
@> \Phi >> R=H_{2*}(\{pt\},R)\:,
\end{CD} \end{equation}
with $\Omega^{U}_{*}$ the {\em cobordism ring of stable almost complex manifolds\/}.
But $\Omega^{U}_{*}\otimes\bb{Q}$ is a polynomial ring in the classes
of all complex projective spaces. Moreover, the characteristic class transformation $cl^{*}$ is uniquely fixed by 
\[\Phi([\bb{P}^{n}(\bb{C})]):= \Phi([\bb{P}^{n}(\bb{C})\to \{pt\}])
= \int_{[P^{n}(C)]} \: (cl^{*}(T\bb{P}^{n}(\bb{C}))\cap [\bb{P}^{n}(\bb{C})])  \] 
for all $n$.
But if $\Phi$ also factorizes over $K_{0}(var/\{pt\})$ then we get by ``additivity''
and ``multiplicativity'':

\begin{equation} \label{eq:genus2}
\Phi([\bb{P}^{n}(\bb{C})]) = 1+ (-y) + \cdots + (-y)^{n} 
\quad \text{with} \quad y:= 1- \Phi([\bb{P}^{1}(\bb{C})]) \:.
\end{equation}
So $\Phi$ is a specialisation of the {\em Hirzebruch $\chi_{y}$-genus\/}
corresponding to the {\em Hirzebruch characteristic class\/} $T^{*}_{y}$ \cite{Hi},
and our {\em Hirzebruch class transformation\/} $T_{y*}$ is the most general 
``additive'' one for homology with values in a $\bb{Q}$-algebra $R$ !\\

Note that the specialisation $y=1$ corresponding to the {\em signature
genus\/} $sign=\chi_{1}$ and the {\em characteristic $L$-class transformation\/} $cl^{*}=L^{*}=T^{*}_{1}$
is the only one that factorizes by the canonical map
$\Omega^{U}_{*}\to \Omega^{SO}_{*}$ over the {\em cobordism ring
$\Omega^{SO}_{*}$ of oriented manifolds\/}:
\begin{displaymath}  \begin{CD} 
Iso^{pro}(sm/\{pt\}) @>>> \Omega^{SO}_{*}\otimes\bb{Q}  \\ 
@VVV  @VV \Phi V\\
K_{0}(var/\{pt\}) 
@> \Phi >> R=H_{2*}(\{pt\},R)\:,
\end{CD} \end{displaymath}
since $[\bb{P}^{1}(\bb{C})]= 0 \in \Omega^{SO}_{*}$ ! \\

More generally, $\Omega^{U}_{*}$ is a commutative ring generated by the classes of all complex projective spaces together with the classes of the {\em Milnor manifolds\/}
\[H_{n,m}:=\{([x_{0},\dots,x_{n}],[y_{0},\dots,y_{m}]) \in 
\bb{P}^{n}(\bb{C}) \times \bb{P}^{m}(\bb{C})|\: \sum \nolimits_{i=0}^{m}\: 
x_{i}\cdot y_{i}= 0 \} \]
for $0<m \leq n$.
These are algebraic (Zariski trivial) projective bundles over projective spaces.
So again any ring homomorphism $\Phi: Iso^{pro}(sm/\{pt\}) \to R$
to a commutative ring $R$, which factorizes as
\begin{displaymath}  \begin{CD} 
Iso^{pro}(sm/\{pt\}) @>>> \Omega^{U}_{*}  \\ 
@VVV  @VV \Phi V\\
K_{0}(var/\{pt\}) @> \Phi >> R
\end{CD} \end{displaymath}
has to be a specialisation of the Hirzebruch $\chi_{y}$-genus.\\
 
If one considers more general (co)homology theories $H^{*}$, it becomes more natural to weaken the {\em normalization condition\/} to 
\begin{itemize}
\item $cl^{*}(L)=f(c^{1}(L))$ for any line bundle $L$, with
$f(z)\in H^{*}(\{pt\})[[z]]$ a formal power series with $f(0)\in H^{*}(\{pt\})$
a {\em unit\/}, 
\end{itemize}
with $cl^{*}$ a multiplicative transformation on the set of isomorphism classes of vector bundles. Note that this is indeed the case for 
\[cl^{*}(V):= \sum_{i=0}^{rk(V)} \; [\Lambda^{i} V^{\vee}] \cdot y^{i}
= \lambda_{y}([V^{\vee}]) \in K^{0}(X)\otimes\bb{Z}[y,(1+y)^{-1}] \]
(with $V^{\vee}$ the dual bundle)
corresponding to our {\em motivic Chern class transformation\/} $mC_{*}$
(with $c^{1}(L)= 1- [L^{\vee}]$ for a line bundle $L$ so that $f(z)= 1+y - y\cdot z$).
But for the composed transformation
\[mC_{*}: K_{0}(var/X) \to G_{0}(X)\otimes\bb{Z}[y]
\to G_{0}(X)\otimes\bb{Z}[y,(1+y)^{-1}] \]
we are not allowed anymore to specialize to the special value $y=-1$ !
And indeed the corresponding genus $\chi_{-1}$ is just the {\em topological
Euler characteristic\/} corresponding to the total Chern class transformation
$cl^{*}=c^{*}$ related to the {\em Chern-Schwartz-MacPherson class transformation\/}
\[c_{*}: K_{0}(var/X) \to A_{*}(X) \to H^{BM}_{2*}(X,\bb{Z}) \:.\]

Here one can deduce this Chern class transformation on the relative
Grothen\-dieck group $K_{0}(var/X)$ without appealing to MacPherson's theorem,
since the distinguished element 
\[d_{X}:= c^{*}(TX)\cap [X] \in A_{*}(X)\]
of a smooth space $X$ satisfies the assumptions of Corollary \ref{cor:bl}.
Condition (1) follows from the projection formula, and condition (2) is an
easy application (by pushing down to $X$) of the classical ``blowing up formula
for Chern classes'' \cite[thm.15.4]{Ful}. In fact this formula is true over any
base field $k$. In particular for $k=\bb{R}$ we get a commutative diagram
of natural transformations:

\begin{equation} \label{eq:SW} \begin{CD} 
K_{0}(var/X) @> e_{2} >> F(X(\bb{R}),\bb{Z}_{2})   \\ 
@VV c_{*} V  @VV w_{*} V\\
A_{*}(X/\bb{R}) @> cl_{R} >> H^{BM}_{*}(X(\bb{R}),\bb{Z}_{2}) \:,
\end{CD} \end{equation}
with $cl_{R}$ the corresponding {\em fundamental class map\/} of Borel-Haefliger \cite{BH},
$F(X(\bb{R}),\bb{Z}_{2})$ the group of real algebraically $\bb{Z}_{2}$-valued
constructible functions and $w_{*}$ the {\em Stiefel-Whitney homology class transformation\/} 
of Sullivan \cite{Su,FuM}. Here $e_{2}$ is again given by
\[e_{2}([f: Y\to X]):=f_{!}1_{Y(R)} \in F(X(\bb{R}),\bb{Z}_{2}) \:,\]
i.e. by taking fiberwise the topological mod 2 Euler characteristic with compact support of the corresponding map $f: Y(\bb{R})\to X(\bb{R})$ on the set of
real points. We also have a canonical fundamental class map
\[Iso^{pro}(sm/\{spec(\bb{R})\}) \to  \Omega^{O}_{*} \:;
[X\to \{pt\}]\mapsto [X(\bb{R})] \]
to the {\em cobordism ring $\Omega^{O}_{*}$ of unoriented manifolds\/}, which is a polynomial
$\bb{Z}_{2}$-algebra generated by the images of all projective spaces
and Milnor manifolds (which are defined over $spec(\bb{R})$). 
So again any ring homomorphism 
\[\Phi: Iso^{pro}(sm/\{spec(\bb{R})\}) \to \bb{Z}_{2}\]
which factorizes as
\begin{displaymath}  \begin{CD} 
Iso^{pro}(sm/\{spec(\bb{R})\}) @>>> \Omega^{O}_{*}  \\ 
@VVV  @VV \Phi V\\
K_{0}(var/\{pt\}) @> \Phi >> \bb{Z}_{2}
\end{CD} \end{displaymath}
has to be the spezialisation $\chi_{1}$ mod 2 of the Hirzebruch $\chi_{y}$-genus,
since now $[\bb{P}^{1}(\bb{R})]= 0 \in \Omega^{O}_{*}$ !

And this is just the mod 2 {\em Euler characteristic\/} corresponding to 
the total Stiefel-Whitney class transformation
$cl^{*}=w^{*}$ for real vector bundles related to 
Sullivan's homology class transformation $w_{*}$.\\

In this way the relative Grothendieck group $K_{0}(var/X)$ of algebraic spaces
unifies all known functorial homology class transformations
\[mC_{*}, T_{y*}, c_{*}, td_{*}, L_{*} \quad \text{and} \quad w_{*} \:,\]
and the ``additivity'' condition ``singles these out'' of all multiplicative
(and normalized) cohomology class transformations $cl^{*}$ for vector bundles.
More precisely, this fixes the {\em genus\/} $\Phi$, but in general not the corresponding
cohomology class $cl^{*}$. This is related to suitable {\em (co)homology operations\/}
on $K^{0}(\cdot)\otimes R$ or $A_{*}(\cdot)\otimes R$, e.g.:
\begin{itemize}
\item The {\em duality involution\/} $[V]\mapsto [V^{\vee}]$ on $K^{0}(\cdot)$,
or the {\em Adams operation\/} $\Psi^{j}$ on $K^{0}(\cdot)\otimes \bb{Z}[1/j]$
(for $j\in \bb{Z}$, and compare with \cite{FL}).
\item The {\em Steenrod $p$-th power operation\/} on $A_{*}(\cdot)\otimes \bb{F}_{p}$
(for $p\in \bb{N}$ a prime number, and compare with \cite{Bro}). 
\end{itemize}

\begin{rem}
The ``additivity'' of the {\em Euler characteristic\/} (with compact support) is well known.
Compare with  \cite[Ex.15.2.10]{Ful} for a characterization
of  the {\em arithmetic genus\/} in terms of a different ``additivity''
property, and with \cite{J} for a characterization of the {\em signature\/}
in terms of ``Novikov additivity'' (as developed in \cite{Si} in the context of
``Witt spaces'').
\end{rem}  
 
Let us finally remark that most of our constructions would apply to a perfect
base field $k$ of positive characteristic, once a suitable version of
resolution of singularities is available (e.g. such that the construction of
the Du Bois complex from \cite{GNA} applies).\\ 

It is a pleasure to thank P.Aluffi for some discussions
on this subject. The paper \cite{A1} was a strong motivation for our work,
which started with the papers \cite{Y1}-\cite{Y6} of the third author.
Some of these papers are partly motivated by the final remark of
MacPherson's survey article:
\begin{itemize} 
\item \cite[p.326]{M2}: ``... It remains to be seen wether there is a unified
theory of characteristic classes of singular varieties like the classical one
outlined above. ...''.
\end{itemize}
We hope that our results give some key to MacPherson's
question and answer  the following question or problem:
\begin{itemize} 
\item \cite[p.267]{Y4}: ``... Is there a theory of characteristic homology classes
unifying the above three characteristic homology classes of possibly
singular varieties? ...''.
\item \cite[p.3367]{A1}: ``... There is a strong motivic feel to the theory of
Chern-Schwartz-Mac\-Pherson classes, although this does not seem to have yet
been congealed into a precise statement in the literature. ...''.
\end{itemize}

\tableofcontents

\section{The generalized Hirzebruch theorem} \label{sec:gHHR}

First we recall the classical {\em generalized Hirzebruch Riemann-Roch theorem\/}
\cite{Hi} (compare with \cite{Y3,Y4}). Let $X$ be a smooth complex projective
variety and $E$ a holomorphic vector bundle over $X$. The
{\em $\chi_{y}$-characteristic\/} of $E$ is defined by \begin{align*}
\chi_{y}(X,E):= 
&\sum_{p\geq 0} \chi(X,E\otimes \Lambda^{p}T^{*}X)\cdot y^{p}\\
=&\sum_{p\geq 0} \left( \sum_{i\geq 0}
(-1)^{i}dim_{C}H^{i}(X,E\otimes \Lambda^{p}T^{*}X) \right)\cdot y^{p} \:,
\end{align*}
with $T^*X$ the holomorphic cotangent bundle of $X$. Then one gets
\begin{equation} \label{eq:gHRR} 
\chi_{y}(X,E)= \int_{X} T^{*}_{y}(TX)\cdot ch_{(1+y)}(E) \cap [X]
\quad \in \bb{Q}[y], \tag{gHRR}
\end{equation}
\[\text{with} \quad ch_{(1+y)}(E):= \sum_{j=1}^{rk\;E} e^{\beta_{j}(1+y)}
\quad \text{and} \quad T^{*}_{y}(TX):= \prod_{i=1}^{dim X}
Q_{y}(\alpha_{i})  \:.\]
Here $\beta_{j}$ are the Chern roots of $E$ and $\alpha_{i}$ are the Chern
roots of the tangent bundle $TX$. Finally $Q_{y}(\alpha)$ is the normalized
power series 
\[Q_{y}(\alpha):= \frac{\alpha(1+y)}{1-e^{-\alpha(1+y)}} -\alpha y
\quad \in \bb{Q}[y][[\alpha]] \:. \]
So this power series $Q_{y}(\alpha)$ specializes to
\begin{displaymath}
Q_{y}(\alpha) = 
\begin{cases}
\:1+\alpha &\text{for $y=-1$,}\\
\:\frac{\alpha}{1-e^{-\alpha}} &\text{for $y=0$,}\\
\:\frac{\alpha}{\tanh \alpha} &\text{for $y=1$.}
\end{cases} \end{displaymath}

Therefore the {\em modified Todd class\/} $T^{*}_{y}(TX)$ unifies the
following important characteristic cohomology classes of $TX$:
\begin{displaymath}
T^{*}_{y}(TX) = 
\begin{cases}
\:c^{*}(TX) &\text{the total {\em Chern class\/} for $y=-1$,}\\
\:td^{*}(TX) &\text{the total {\em Todd class\/} for $y=0$,}\\
\:L^{*}(TX) &\text{the total {\em Thom-Hirzebruch L-class\/} for $y=1$.}
\end{cases} \end{displaymath}

Note that (\ref{eq:gHRR}) implies for $y=0$ the classical {\em Hirzebruch 
Riemann-Roch theorem\/} \cite{Hi}:
\begin{equation} \label{eq:HHR} 
\chi(X,E)= \int_{X} td^{*}(TX)\cdot ch^{*}(E) \cap [X] \:, \tag{HRR}
\end{equation}
\[\text{with} \quad ch^{*}(E):= \sum_{j=1}^{rk\;E} e^{\beta_{j}}\]
the classical {\em Chern character\/}. This is a ring homomorphism (for the tensor product)
\[ch^{*}: K^{0}(X) \to H^{2*}(X)\:,\]
so that one gets back (\ref{eq:gHRR}) by
\begin{equation} \label{eq:gHRR'} \tag{gHRR'}
\begin{split} 
\chi_{y}(X,E) &= \int_{X} td^{*}(TX)\cdot ch^{*}(E\otimes \lambda_{y}(T^{*}X)) \cap [X] \\
 &= \int_{X} td^{*}(TX)\cdot ch^{*}(\lambda_{y}(T^{*}X))\cdot ch^{*}(E)  \cap [X] \:,
\end{split} 
\end{equation}
with $\widetilde{T}^{*}_{y}(TX):= td^{*}(TX)\cdot ch^{*}(\lambda_{y}(T^{*}X))$
corresponding to the unnormalized power series (compare \cite[p.11,p.61]{HBJ})
\[\widetilde{Q}_{y}(\alpha):= \frac{\alpha(1+ye^{-\alpha})}{1-e^{-\alpha}} 
\quad \in \bb{Q}[y][[\alpha]] \:, \quad \text{with}\quad
\widetilde{Q}_{y}(0)= 1+y \:. \]

Then the relation (compare \cite[p.62]{HBJ})
\[Q_{y}(\alpha)=\widetilde{Q}_{y}(\alpha\cdot (1+y)) \cdot (1+y)^{-1} \in \bb{Q}[y][[\alpha]] \]
implies for $X$ pure $d$-dimensional:
\begin{equation} \label{eq:comp}
\begin{split}
T_{y}^{i}(TX) = (1+y)^{i-d} \cdot \widetilde{T}^{i}_{y}(TX) \in H^{2i}(X)\otimes \bb{Q}[y]
\quad \text{and} \\
(T_{y})_{j}(X) = (1+y)^{-j} \cdot (\widetilde{T}_{y})_{j}(X) \in H_{2j}(X)\otimes \bb{Q}[y] \:.
\end{split} \end{equation}
Here we use the notation $cl_{*}(X):=cl^{*}(TX)\cap [X]$ for a characteristic class $cl^{*}$
of vector bundles. In particular $(T_{y})_{*}(X)$ and $(\widetilde{T}_{y})_{*}(X)$ agree in degree $0$ 
so that
\[\chi_{y}(X):=\chi_{y}(X,{\cal O}_{X})= \int_{X} \: (T_{y})_{*}(X)
= \int_{X} \: (\widetilde{T}_{y})_{*}(X) \:.\]
So this ``twisting'' by powers of $1+y$ just comes from changing $\widetilde{Q}_{y}$
to the normalized power series $Q_{y}$, which has the right specialization properties.
 
And (\ref{eq:gHRR'}) implies (\ref{eq:gHRR}) by a similar calculation.\\

As an example, Hirzebruch \cite[lem.1.8.1]{Hi} gets by a simple residue calculation the
equation
\[ \int_{X} \: (T_{y})_{*}(X) = 1+ (-y) + \cdots + (-y)^{n} \quad \text{for}
\quad X=\bb{P}^{n}(\bb{C}) \:.\]
The gRRH-theorem for the trivial bundle $E$ specializes to
the calculation of the following important invariants:
\begin{displaymath}
\chi_{y}(X) = 
\begin{cases}
\: e(X) = \int_{X} c^{*}(TX)\cap [X] &\text{the {\em Euler
characteristic\/} for $y=-1$,}\\ 
\: \chi(X) = \int_{X} td^{*}(TX)\cap [X] &\text{the {\em arithmetic genus\/} for
$y=0$,}\\ 
\: sign(X) = \int_{X} L^{*}(TX)\cap [X]  &\text{the {\em signature\/} for $y=1$,}
\end{cases} \end{displaymath}
corresponding to the {\em Poincar\'e-Hopf or Gauss-Bonnet theorem\/} ($y=-1$), the
{\em Hirzebruch Riemann-Roch theorem\/} ($y=0$) and the
{\em Hirzebruch signature theorem\/} ($y=1$).\\

These three invariants and classes have been generalized to a singular complex
algebraic variety $X$ in the following way (where the invariants are only
defined for $X$ compact):

\begin{equation} \label{eq:y-1} \tag{$y=-1$}
e(X) = \int_{X} c_{*}(X), \:\: \text{with} \:\:
c_{*}: F(X) \to H_{*}(X):=
\begin{cases}
A_{*}(X)\\
H^{BM}_{2*}(X,\bb{Z})
\end{cases}
\end{equation}
the {\em Chern class transformation\/} of MacPherson \cite{M1,Ken}
from the abelian group $F(X)$
of complex algebraically constructible functions to homology,
where one can use Chow groups $A_{*}(\cdot)$ or Borel-Moore homology groups
$H^{BM}_{2*}(\cdot,\bb{Z})$. Then $c_{*}(X):=c_{*}(1_{X})$
agrees by \cite{BrS} via ``Alexander duality'' for compact $X$ embeddable into a complex manifold with
the  {\em Schwartz class\/} of $X$  as introduced before by
M.-H. Schwartz \cite{Schwa}.

\begin{equation} \label{eq:y0} \tag{$y=0$}
\chi(X) = \int_{X} td_{*}(X), \quad \text{with} \quad
td_{*}: G_{0}(X) \to H_{*}(X)\otimes \bb{Q} 
\end{equation}
the {\em Todd transformation\/} in the singular Riemann-Roch theorem of
Baum-Fulton-MacPherson \cite{BFM,FM} (for Borel-Moore homology) or Fulton \cite{Ful}
(for Chow groups). Here $G_{0}(X)$ is the Grothendieck group
of coherent sheaves.
Then $td_{*}(X):=td_{*}([{\cal O}_{X}])$,
with $[{\cal O}_{X}]$ the class of the structure sheaf.\\
 
Finally for compact $X$ one also has
\begin{equation} \label{eq:y1} \tag{$y=1$}
sign(X) = \int_{X} L_{*}(X), \quad \text{with} \quad
L_{*}: \Omega(X) \to H_{2*}(X,\bb{Q}) 
\end{equation}
the {\em homology L-class transformation\/} of Cappell-Shaneson \cite{CS1}
(as reformulated by Yokura \cite{Y1} and corrected in Section \ref{sec:Lclass}). 
Here $\Omega(X)$ is the abelian group
of cobordism classes of selfdual constructible complexes.
Then $L_{*}(X):=L_{*}([{\cal IC}_{X}])$ is the {\em homology L-class\/} of
Goresky-MacPherson \cite{GM},
with $[{\cal IC}_{X}]$ the class of their intersection cohomology complex.
For a {\em rational PL-homology manifold\/} $X$,
these L-classes are due to Thom \cite{Thom} (compare \cite[sec.20]{MS}).

\begin{rem} \label{rem:AS}
The discussion above applies to any compact complex manifold $X$,
since the generalized Hirzebruch Riemann-Roch theorem is also true
in this context by an application of the {\em Atiyah-Singer Index theorem\/}
\cite{AS}. Similarly, (\ref{eq:gHRR}) follows as before in the algebraic context over
any base field $k$ from the corresponding {\em HRR-theorem\/} \cite[cor.15.2.1]{Ful}
for Chow-groups $A_{*}(\cdot)\otimes \bb{Q}$
(instead of homology $H_{2*}(\cdot,\bb{Q})$). And this HRR-theorem is just
the special case of the {\em Grothendieck Riemann-Roch theorem\/} \cite[thm.15.2]{Ful}
for a constant map $X\to spec(k)$, with $X$ a smooth complete variety.
\end{rem}

All these transformations commute with the corresponding pushdown
for proper maps (where all spaces are assumed to be compact in the case
of the L-class transformation). They are uniquely characterized
by this pushdown property and the {\em normalization condition\/} that for $X$ smooth
and pure-dimensional one gets back the corresponding classes of $TX$:
\[c_{*}(X)=c^{*}(TX)\cap [X],\quad td_{*}(X)=Td^{*}(TX)\cap [X] \:\:
\text{and} \:\: L_{*}(X)=L^{*}(TX)\cap [X] \:.\]
Here the uniqueness result follows from resolution of singularities,
and in the case of the $L$-class transformation one has to be more careful:
This normalization fixes $L_{*}$ only on the image of the transformation
$sd$ from (\ref{eq:real2})!\\

So all these theories have the same formalism, but they are defined 
on completely different theories! Nevertheless, it is natural to ask
for another theory of characteristic homology classes of 
singular complex algebraic varieties, which unifies the above characteristic
homology class transformations (as in \cite{M2,Y3,Y4}). 
Of course in the smooth case, this is done by the 
{\em generalized Todd class\/} $T^{*}_{y}(TX) \cap [X]$ of the tangent
bundle. We now explain our solution to this question.

\section{Motivic Chern classes for singular varieties}

In the following we consider reduced seperated schemes of finite type
over a base field $k$ of characteristic $0$, and for simplicity we just call them algebraic varieties.  \\

Let $K_{0}(var/X)$ be the {\em relative Grothendieck group of algebraic
varieties over\/} $X$, i.e. the quotient of the free abelian group of isomorphism
classes of algebraic morphisms $Y\to X$ to $X$,  modulo the ``additivity
relation'' (\ref{eq:add}).
These relative groups were introduced by Looijenga
in his Bourbaki talk \cite{Lo} about motivic measures and motivic integration,
and then further studied by Bittner \cite{Bi}. From our point of view,
these are a ``motivic version'' of the group $F(X)$ of algebraically constructible
functions (compare also with \cite{CL}). 

In particular, they have the same formalism, i.e.
functorial {\em pushdown\/} $f_{!}$ and {\em pullback\/} $f^{*}$ for any algebraic map
$f:X'\to X$ (which is not necessarily proper), together with a ring
multiplication (with unit $[id_{X}]=k^{*}[id_{pt}]$ for $k:X\to \{pt\}$ a constant
map) satisfying the {\em projection formula\/} 
\[f_{!}(\alpha\cdot f^{*}\beta)=(f_{!}\alpha)\cdot \beta\]
and the {\em base change formula\/} $g^{*}f_{!}=f'_{!}g'^{*}$ for any cartesian
diagram
\begin{displaymath} \begin{CD}
Y' @>g'>> X' \\
@V f' VV  @VV f V \\
Y @> g >> X \:. 
\end{CD} \end{displaymath}
For later use, let us recall the simple definition of the {\em pullback\/} and
{\em pushdown\/}  for $f: X'\to X$, and of {\em exterior products\/}:
\[f_{!}([h:Z\to X'])=[f\circ h: Z\to X] \quad \text{and}\]
\[[Z\to X]\times [Z'\to X'] = [Z\times Z' \to X \times X']\:.\]
\[\text{Moreover} \quad f^{*}([g:Y\to X])= [g': Y'\to X']\]
is defined by taking fiber products as above.\\

By these exterior products, $K_{0}(var/\{pt\})$ becomes a commutative ring and
$K_{0}(var/X)$ a $K_{0}(var/\{pt\})$-module such that $f_{!}$ and $f^{*}$ are
$K_{0}(var/\{pt\})$-linear. Moreover, the homomorphism $e$ from (\ref{eq:e})
commutes with all these transformation. In particular it is a ring
homomorphism.

\begin{rem} \label{rem:iso}
If we consider only {\em proper\/} morphism for the pushdown,
and only {\em smooth\/} morphism (of constant relative dimension) for the pullback,
then the same formalism applies to the groups $Iso^{pr}(var/X)$ and
$Iso^{pr}(sm/X)$. But here we only have an {\em exterior product\/},
but no ring structure on these groups, since the diagonal morphism
$X\to X\times X$ is not a smooth morphism. But this is true for $X=\{pt\}=spec(k)$
a point, so $Iso^{pr}(var/\{pt\})$ and $Iso^{pr}(sm/\{pt\})$ become a commutative ring
such that the transformations above are linear over these rings.
Moreover, the group epimorphisms
\[Iso^{pr}(sm/X) \to Iso^{pr}(var/X) \to K_{0}(var/X)\]
commute with these three operations. A similar remark also applies
to the compactifiable complex analytic context.
\end{rem}

Note that the groups $Iso^{pr}(var/X)$ and $Iso^{pr}(sm/X)$ are {\em graded\/} by the 
dimension of the spaces mapping to $X$
(if we work with pure dimensional spaces), whereas $K_{0}(var/X)$ only becomes
a {\em filtered group\/} with $F_{k}K_{0}(var/X)$ generated by $[X'\to X]$ with
$dim X' \leq k$.\\

Let us now explain the simple 
\begin{proof}[Proof of lemma \ref{lem:iso1}.]
The surjective projection $\pi: Iso^{pr}(var/X) \to K_{0}(var/X)$ factorizes by definition over the
quotient $Iso^{pr}(var/X)/(ac)$. So it is enough to define an inverse map
\[\phi: K_{0}(var/X) \to  Iso^{pr}(var/X)/(ac)\:.\]
And this is induced by the map 
\[\begin{split}
\phi: Iso(var/X) \to  Iso^{pr}(var/X)/(ac)\:; \\
[X'\to X] \mapsto [\bar{X}'\to X]_{ac} - [\bar{X}'\backslash X'\to X]_{ac}
\end{split} \]
taking the ``difference class'' for a compactification $X'\subset \bar{X}'\to X$ of the corresponding map,
with $X'\subset \bar{X}'$ (Zariski) open and $\bar{X}'\to X$ proper (i.e. by Nagata's theorem \cite{L}).
By using the fiber product of two such compactifications, one gets from
the relation (\ref{eq:ac}) that this is well defined. A simple calculation shows that it satisfies also the ``additivity relation'' (\ref{eq:add}) for $Z\subset X'$ closed:
\[\begin{split}
\phi\bigl([X'\to X]\bigl) = &[\bar{X}'\to X]_{ac} - [\bar{X}'\backslash X'\to X]_{ac} \\
= [\bar{X}'\to X]_{ac} - &[Z\cup \bar{X}'\backslash X'\to X]_{ac} 
+ [Z\cup \bar{X}'\backslash X'\to X]_{ac} - [\bar{X}'\backslash X'\to X]_{ac}\\
= &\phi\bigl([X'\backslash Z\to X]\bigl) + \phi\bigl([Z\to X]\bigl) \:.
\end{split} \]
For $X'\to X$ proper we can take $X'=\bar{X}'$ so that $\phi\circ \pi = id$.

One can avoid the use of Nagata's theorem by decomposing $X'$ first into the disjoint
union of quasi-projective (e.g affine) pieces $X'_{i}$. For the $X'_{i}$ one uses the definition as before (with $\bar{X}'\to X$ projective) and then one defines
\[\phi\bigl([X'\to X]\bigl) := \sum\nolimits_{i}\: \phi\bigl([X_{i}'\to X]\bigl) \:,\]
which is well defined by ``additivity''.
\end{proof}

\begin{rem} The proof of Theorem \ref{thm:blowup} given in \cite{Bi} is of similar nature.
But here one takes for $X'$ smooth a compactification such that $\bar{X}'\backslash X'$
is a divisor with simple normal crossings. This can be done by Hironaka's
resolution of singularities. Then one has to compare two such compactifications.
And for this the ``weak factorization theorem'' of \cite{AKMW,W} is used,
which relates two such compactifications by a finite sequence of blowing ups and
blowing downs along (suitable) smooth centers! Moreover, it is enough to consider
the ``blow-up'' relation only for
(pure dimensional) quasi-projective varieties $X'$ with $X'\to X$ a projective
morphism.
\end{rem}

By the discussion in the introduction we therefore get the

\begin{thm} \label{thm:mC}
There exists a unique group homomorphism $mC_{*}$ commuting with pushdown
for proper maps:
\[mC_{*}: K_{0}(var/X)\to G_{0}(X)\otimes\bb{Z}[y] \:,\]
satisfying the normalization condition
\[mC_{*}([id_{X}])= \sum_{i=0}^{dim X} \; [\Lambda^{i} T^{*}X]\cdot y^{i}
= \lambda_{y}([T^{*}X])\cap [{\cal O}_{X}]\]
for $X$ smooth and pure-dimensional.
Here $\lambda_{y}: K^{0}(X)\to K^{0}(X)\otimes\bb{Z}[[y]]$ is the {\em total
$\lambda$-class transformation\/} on the Grothendieck 
$K^{0}(X)$ of coherent locally free sheaves on $X$, with
$\cap [{\cal O}_{X}]: K^{0}(X)\to G_{0}(X)$ induced by 
$\otimes {\cal O}_{X}$, which is an isomorphism for $X$ smooth.
\hfill $\Box$
\end{thm}

\begin{rem} For a compact complex space $X$ this result is true for $K_{0}(an/X)$,
except that $\cap [{\cal O}_{X}]: K^{0}(X)\to G_{0}(X)$ needs not to be an isomorphism.
\end{rem}

\begin{cor} \label{cor:mC}
\begin{enumerate}
\item $mC_{*}$ is filtration preserving, if $G_{0}(X)\otimes\bb{Z}[y]$
has the induced filtration coming from the grading with $y$ of degree one.
\item $mC_{0}: K_{0}(var/X)\to G_{0}(X) $ is the unique group homomorphism
 commuting with pushdown for proper maps and 
satisfying the normalization condition
$mC_{0}([id_{X}])=[{\cal O}_{X}]$ 
for $X$ smooth (and pure-dimensional).
\item $mC_{*}$ commutes with exterior products.
\item One has the following {\em Verdier Riemann-Roch formula\/} for $f:X'\to X$
a smooth morphism (of constant relative dimension):
\[\lambda_{y}(T^{*}_{f}) \cap f^{*}mC_{*}([Z\to X]) 
= mC_{*}f^{*}([Z\to X]) \:.\]
Here $T_{f}$ is the bundle on $X'$ of tangent spaces to the fibers of $f$,
i.e. the kernel of the surjection $df: TX'\to f^{*}TX$.
Moreover $f^{*}: G_{0}(X)\otimes\bb{Z}[y] \to
G_{0}(X')\otimes\bb{Z}[y]$
is induced from the corresponding pullback of Grothendieck groups
by linear extension over $\bb{Z}[y]$.
In particular $mC_{*}$ commutes with pullback under \'{e}tale morphisms
(i.e. smooth of relative dimension $0$).
\end{enumerate} \end{cor}

\begin{proof} \begin{enumerate}
\item follows by induction on $dim\: X$ from resolution of
singularities, ``additivity'' and the normalization condition.
\item follows from the fact that pushdown on $G_{0}\otimes\bb{Z}[y]$
is degree preserving. 
\item follows  from the  normalization condition together with
\[(f\times f')_{*} \simeq f_{*}\times f'_{*}
\quad \text{on $K_{0}(var/X\times X')$ and $G_{0}(X\times X')
\otimes\bb{Z}[y]$}\] 
for $f,f'$ proper, and by the multiplicativity of
$\lambda_{y}((\cdot)^{\vee})$: 
\[\lambda_{y}(T^{*}(X\times X')) = 
\lambda_{y}(T^{*}X) \times \lambda_{y}(T^{*}X') \:.\]
Compare also with \cite{Kw,KwY} for the case of Chern classes. 
\item It is enough to prove the claim for $g: Z\to X$ proper
with $Z$ smooth (and pure dimensional). Then it
follows from the   projection formula 
\[g'_{*}(\alpha\cdot g'^{*}\beta)=(g'_{*}\alpha)\cdot \beta\] 
for $g'$ proper and the base
change formula $f^{*}g_{*}=g'_{*}f'^{*}$ for the cartesian diagram
\begin{displaymath} \begin{CD}
Z' @>f'>> Z \\
@V g' VV  @VV g V \\
X' @> f >> X 
\end{CD} \end{displaymath}
with $g,g'$ proper and $f,f'$ smooth (of constant relative dimension).
Here these formulae also hold for $G_{0}(\cdot)\otimes\bb{Z}[y]$:
\begin{align*}
&\lambda_{y}(T^{*}_{f}) \cap f^{*}mC_{*}([Z\to X]) =
\lambda_{y}(T^{*}_{f}) \cap f^{*}g_{*}mC_{*}([id_{Z}]) \\ 
=&\lambda_{y}(T^{*}_{f}) \cap g'_{*}f'^{*}mC_{*}([id_{Z}]) 
= g'_{*}\left( g'^{*}\lambda_{y}(T^{*}_{f})\cap f'^{*}(
\lambda_{y}(T^{*}Z)\cap [{\cal O}_{Z}])\right)\\
=& g'_{*}\left( \lambda_{y}(T^{*}_{f'})\cup f'^{*}
\lambda_{y}(T^{*}Z)\cap [{\cal O}_{Z'}]\right)
= g'_{*} \left( \lambda_{y}(T^{*}Z')\cap [{\cal O}_{Z'}] \right)\\
=& g'_{*}mC_{*}([id_{Z'}]) = mC_{*}([Z'\to X']) =  
mC_{*}f^{*}([Z\to X]) \:.
\end{align*} 
Of course, we also used $f'^{*}[{\cal O}_{Z}]=[{\cal O}_{Z'}]$
together with the multiplicativity and functoriality
of $\lambda_{y}((\cdot)^{\vee})$.
Compare also with \cite{Y5} for the case of Chern classes.
\end{enumerate}
\end{proof}

\section{Hirzebruch classes for singular varieties}

We continue to work with reduced algebraic varieties over a base field $k$
of characteristic $0$. In the following, the homology group 
$H_{*}(\cdot)$ is either the Chow group or the Borel-Moore homology group
(in even degrees) in case $k=\bb{C}$.\\

Let us start with a reformulation by Yokura \cite{Y4} of the {\em singular Riemann-Roch theorem\/} of Baum-Fulton-MacPherson (for Borel-Moore homology)
and Fulton (for Chow groups).
 Consider the group homomorphism
\begin{equation} \label{eq:gBFM} \tag{gBFM}
\begin{split}
td_{(1+y)}: G_{0}(X)\otimes\bb{Z}[y] &\to H_{*}(X)\otimes\bb{Q}[y,(1+y)^{-1}]
\:,\\
td_{(1+y)}([{\cal F}])&:= \sum_{i\geq 0}\: td_{i}([{\cal F}])\cdot (1+y)^{-i}
\:, \end{split} \end{equation}
with $td_{i}$ the degree $i$ component of the transformation $td_{*}$,
which is linearly extended over $\bb{Z}[y]$. Since $td_{*}$ is
degree preserving, this new transformation also commutes with proper pushdown
(which again is defined by linear extension over $\bb{Z}[y]$).
By \cite{Y3} we have the

\begin{lem} \label{lem:Yok}
Assume $X$ is smooth (and pure dimensional). Then
\[td_{(1+y)}\left( \lambda_{y}(T^{*}X) \cap [{\cal O}_{X}] \right) = T_{y}^{*}(TX)
\cap [X] \in H_{*}(X)\otimes\bb{Q}[y]\:,\]
\[td_{*}\left( \lambda_{y}(T^{*}X) \cap [{\cal O}_{X}] \right) = \tilde{T}_{y}^{*}(TX)
\cap [X] \in H_{*}(X)\otimes\bb{Q}[y] \:.\]
\end{lem}

Let us sketch the simple proof. 
Since $X$ is smooth, the singular Riemann-Roch theorem reduces to the classical
{\em Grothendieck Rie\-mann-Roch theorem\/}, i.e. 
\begin{equation} \label{GRR} \tag{GRR}
\begin{split}
td_{*}(E\cap [{\cal O}_{X}])&= ch^{*}(E)\cup Td^{*}(TX)\cap [{\cal O}_{X}]:\\ G_{0}(X)&=K^{0}(X) \to
H_{*}(X)\otimes \bb{Q}\:,
\end{split} \end{equation}
with $ch^{*}: K^{0}(X) \to H^{*}(X)\otimes \bb{Q}$ the Chern character.
So 
\[td_{*}\left( \lambda_{y}(T^{*}X) \cap [{\cal O}_{X}] \right) = 
(\widetilde{T}_{y})_{*}(X) =  \tilde{T}_{y}^{*}(TX) \cap [X] \]
and the claim follows from (\ref{eq:comp}). In particular, the ``twisting''
(by powers of $1+y$) form $td_{*}$ to $td_{(1+y)}$  is only used to get the right normalization in lemma \ref{lem:Yok}!

\begin{rem} The same reasoning applies also to the {\em compactifiable complex analytic
context\/} with the Todd transformation (\ref{eq:toddan}), since for $X$ smooth the following diagram commutes by \cite{Levy}:
\begin{equation} \begin{CD} \label{eq:Ktop}
K^{0}(X) @> can >> K^{0}_{top}(X)\\
@VV \cap {\cal O}_{X} V @V PD V \cap [X] V \\
G_{0}(X) @>> \alpha > K_{0}^{top}(X) \:.
\end{CD} \end{equation} 
Here $can$ is the canonical map induced from taking the underlying topological
complex vector bundle of a holomorphic vector bundle, and $PD$ is the 
Poincar\'e duality for topological $K$-theory.
\end{rem} 

Now we are ready to introduce our second motivic transformation.
 
\begin{thm} \label{thm:Ty}
There exist unique group homomorphisms $T_{y*},\tilde{T}_{y*}$ commuting with pushdown
for proper maps:
\[T_{y*},\tilde{T}_{y*}: K_{0}(var/X)\to H_{*}(X)\otimes\bb{Q}[y] \:,\]
satisfying the normalization condition
\[T_{y*}([id_{X}])=T^{*}_{y}(TX) \cap [X] \quad , \quad
\tilde{T}_{y*}([id_{X}])=\tilde{T}^{*}_{y}(TX) \cap [X]\]
for $X$ smooth and pure-dimensional.
\end{thm}

\begin{proof} 
The natural transformation $T_{y*}$ is given as the composition
\[T_{y*}:= td_{(1+y)} \circ mC_{*}: K_{0}(var/X)\to H_{*}(X)\otimes\bb{Q}[y]
\subset H_{*}(X)\otimes\bb{Q}[y,(1+y)^{-1}] \:.\]
And similarly
\[\tilde{T}_{y*}:= td_{*} \circ mC_{*}: K_{0}(var/X)\to H_{*}(X)\otimes\bb{Q}[y]
 \:.\]
The normalization condition follows from lemma \ref{lem:Yok}.
Since $Iso^{pr}(sm/X)\to K_{0}(var/X)$ is surjective, we get uniqueness
together with $T_{y*}(K_{0}(var/X))\subset  H_{*}(X)\otimes\bb{Q}[y]$.
\end{proof}

By the same argument one gets from functoriality and normalization the
commutativity of the diagrams (\ref{eq:real1}) and (\ref{eq:real2}).

\begin{rem} The same result and proof also apply to the {\em compactifiable
complex analytic context\/}.
\end{rem}

Before we state some corollaries, let us introduce the following
filtration on $H_{*}(X)\otimes\bb{Q}[y]$:
 \[F_{k}\left(H_{*}(X)\otimes\bb{Q}[y]\right):=
F_{k}H_{*}(X)\otimes F_{k}\bb{Q}[y] \:,\] 
where each factor has its canonical
filtration coming from the natural grading. In particular, any evaluation
homomorphism $H_{*}(X)\otimes\bb{Q}[y]\to H_{*}(X)$ for a $y\in \bb{Q}$
is then also filtration preserving.

\begin{cor} \label{cor:Ty}
\begin{enumerate}
\item $T_{y*}$ is filtration preserving.
\item $T_{y*},\tilde{T}_{y*}$ commute with exterior products.
\item One has the following {\em Verdier Riemann-Roch formula\/} for $f:X'\to X$
a smooth morphism (of constant relative dimension):
\[T^{*}_{y}(T_{f}) \cap f^{*}T_{y*}([Z\to X]) 
= T_{y*}f^{*}([Z\to X]) \:,\]
and similarly for $\tilde{T}_{y*}$.
In particular $T_{y*},\tilde{T}_{y*}$ commute with pullback under \'{e}tale morphisms.
\hfill $\Box$
\end{enumerate}
\end{cor}

Let us recall our definition
\begin{equation} \label{eq:Hclass} \begin{split}
T_{y*}(X)&:= T_{y*}([id_{X}]) \in  H_{*}(X)\otimes\bb{Q}[y]\\
\chi_{y}(X)&:= T_{y}([X\to \{pt\}]) \in \bb{Z}[y]\subset \bb{Q}[y]
\end{split} \end{equation}
for the {\em Hirzebruch class\/} and {\em characteristic\/} of the singular space $X$ (and similarly for the {\em motivic Chern class\/}), so that 
\[T_{-1}([id_{X}]) = c_{*}(X)\otimes \bb{Q} \quad \text{for any singular $X$.}\]

It seems that our Hirzebruch class 
$T_{y*}(X)$ corresponds (for $k=\bb{C}$) to a similar class announced 
some years ago by Cappell and Shaneson \cite{CS2} and \cite[sec.4]{Sh}. 
If this is the case,
then there is a mistake in their announcement, because one of their statements
can be interpreted as claiming that
$T_{0*}([id_{X}])$ is the singular Todd class $td_{*}(X)$ for any singular
$X$. But this is not the case, since 
\[mC_{0}([id_{X}])\neq [{\cal O}_{X}]\in G_{0}(X)
\quad \text{and} \quad T_{0*}([id_{X}]) \neq td_{*}(X) \] 
for some singular spaces $X$.

\begin{ex} \label{ex:T0-Todd}
Let $X$ be a singular curve (i.e. $dim\; X=1$) such that $X$ is not maximal
(maximal is sometimes also called weakly normal), but the weak normalization
$X_{max}$ is smooth. Then the canonical projection $\pi : X':=X_{max}\to X$
is not an isomorphism, but nevertheless a topological homeomorphism.
By ``additivity'' one gets
\[\pi_{*}([id_{X'}]) = [id_{X}] \in K_{0}(var/X) \quad\text{so that}\]
\[T_{0*}(X)=\pi_{*}T_{0*}(X')=\pi_{*}td_{*}(X')=td_{*}\pi_{*}([{\cal O}_{X'}])
\:.\]
But by assumption $\pi_{*}([{\cal O}_{X'}])=[{\cal O}_{X}]+ 
n\cdot [{\cal O}_{pt}]$ with $n>0$ so that
\[T_{0*}(X)= td_{*}(X)+ n\cdot [pt] \neq td_{*}(X) \in H_{*}(X)\otimes \bb{Q}
\:.\]
One gets similar examples in any dimension by taking the product with a
projective space.
\end{ex}

Taking a complete singular curve $X$ over $k=\bb{C}$ such that the
normalization $\pi: X':=X_{nor}\to X$ is not a topological homeomorphism, one
gets in the same way examples of singular $X$ with 
\[T_{1*}(X)\neq L_{*}(X)
:=L_{*}([{\cal IC}_{X}]) \:.\]
Note that the normalization map $\pi$ is a {\em small resolution\/} of singularities
so that $\pi_{*}({\cal IC}_{X'})={\cal IC}_{X}$. To distinguish between
these characteristic classes, we call (for later reasons) $T_{0*}(X)$
the {\em Hodge-Todd class\/} and $T_{1*}(X)$ the {\em Hodge L-class\/} of $X$.

\begin{ex} \label{ex:T0=Todd} 
Assume that the algebraic variety $X$ has
at most ``Du Bois singularities'' in the sense that
\[can: {\cal O}_{X}= gr_{F}^{0}(DR({\cal O}_{X})) \to
gr_{F}^{0}(\underline{\Omega}^{*}_{X})\]
is a quasi-isomorphism.
For example $X$ has only ``rational singularities'' \cite{Kov,Sai5}, e.g. $X$ is a toric variety. Then
\[mC_{0}([id_{X}]) = [{\cal O}_{X}]\in G_{0}(X)\] 
and therefore $T_{0*}(X)=td_{*}(X)$.
\end{ex}

Using the ``additivity'' in $K_{0}(var/X)$ and the natural transformation
$T_{y*}$, one gets similar additivity properties of the Todd class
$td_{*}(X)=T_{0*}(X)$ for $X$ smooth (or with a most ``Du Bois singularities''),
which seem to be new and do not follow directly from the original
definition.

\begin{ex} \label{ex:addTodd}
\begin{enumerate}
\item By the gHRR-theorem we have
$T_{y*}([\bb{P}^{1}\to pt])=\chi_{y}(\bb{P}^{1})=1-y$ so that 
\[\chi_{y}(\bb{A}^{1})=-y \quad \text{and} \quad
\chi_{y}(\bb{P}^{n})= 1 -y + \dots + (-y)^{n} \]
by ``additivity'' and ``multiplicativity'' for exterior products.
\item $T_{y*}$ becomes multiplicative in Zariski locally trivial bundles,
e.g. if $E\to X$ is an algebraic vector bundle of rank $r+1$, then
the corresponding projective bundle $\bb{P}(E)\to X$ is 
Zariski locally trivial so that
\[T_{y*}([\bb{P}(E)\to X])= T_{y*}([id_{X}]) \cdot (1 -y + \dots + (-y)^{r})
\in  H_{*}(X)\otimes \bb{Q}[y] \:.\]
\item Let $\pi: X'\to X$ be the blow-up of an algebraic variety $X$ along
an algebraic subvariety $Y$ such that the inclusion $Y\to X$ is a regular
embedding of pure codimension $r+1$
(e.g. $X$ and $Y$ are smooth). Then $\pi$ is an isomorphism over
$X\backslash Y$ and a projective bundle over $Y$ corresponding to the normal
bundle of $Y$ in $X$ of rank $r+1$. So by ``additivity'' one gets 
\[T_{y*}([X'\to X])= T_{y*}([id_{X}])+ T_{y*}([Y\to X]) \cdot 
(-y + \dots + (-y)^{r}) \]
in  $H_{*}(X)\otimes \bb{Q}[y]$.  
In particular $T_{0*}([X'\to X])= T_{0*}([id_{X}])$, which is a
homology class version of the {\em birational invariance of the
arithmetic genus\/} $\chi_{0}$. More generally, by pushing down to a
point we get for $X$ complete the {\em blow-up formula\/}
\begin{equation} \label{eq:blup}
\chi_{y}(X')= \chi_{y}(X) + \chi_{y}(Y) \cdot 
(-y + \dots + (-y)^{r}) \in  \bb{Q}[y] \:.
\end{equation}
\end{enumerate} \end{ex}

Note that in case 3. one also has $\pi_{*}td_{*}(X')=td_{*}(X)$,
by functoriality of $td_{*}$ and the relation 
$[R\pi_{*}\pi^{*}{\cal O}_{X}]=[{\cal O}_{X}]$ for such a blow-up.
Using the ``weak factorization theorem'' \cite{AKMW,W}, we get the following result,
which seems to be new and was motivated by a corresponding study
of Aluffi about Chern classes \cite{A1}, and of Borisov-Libgober about elliptic classes
\cite{BL2}.

\begin{cor} \label{cor:tdbirat}
Let $\pi: Y\to X$ be a resolution of singularities. Then the class 
\[\pi_{*}(Td^{*}(TY)\cap [Y])=\pi_{*}T_{0}([id_{Y}]) \in H_{*}(X)\otimes
\bb{Q}\]
is independent of $Y$.
\end{cor}

\begin{proof} Let $\pi: Y\to X$ and $\pi': Y'\to X$ be two
resolution of singularities, together with a resolution
of singularities of the fiber-product
$Z\to Y\times _{X} Y'$ so that we get induced birational
morphisms $p: Z\to Y$ and $p': Z\to Y'$.
By the ``weak factorization theorem'' \cite{AKMW,W}, this map $p$ (or $p'$) can be
decomposed as a finite sequence of projections from smooth spaces lying over
$Y$ (or $Y'$), which are obtained by  blowing up or blowing down along  
smooth centers. By the birational invariance above we get
$\pi_{*}td_{*}(Z)= td_{*}(X)$ (or $\pi'_{*}td_{*}(Z)= td_{*}(X')$),
from which the claim follows.
\end{proof}

\begin{rem} \label{rem:local}
In the framework of motivic integration \cite{DL,Lo,Craw,V}, it is natural to localize
the $K_{0}(var/\{pt\})$-module $K_{0}(var/X)$ with respect to the class
of the affine line $[\bb{A}^{1}\to \{pt\}]=: \bb{L}$. Here the module structure
comes from pullback along a constant map $X\to \{pt\}$. Then $mC_{*},T_{y*}$ 
and $\tilde{T}_{y*}$ induce similar transformations on $M(var/X):=K_{0}(var/X)[\bb{L}^{-1}]$:
\[\begin{split}
mC_{*} &: M(var/X)\to G_{0}(X)\otimes \bb{Z}[y,y^{-1}]\:,\\
T_{y*}, \tilde{T}_{y*}
&: M(var/X)\to H_{*}(X)\otimes \bb{Q}[y,y^{-1}]\:,\\
\end{split} \]
since $\chi_{y}(\bb{L})=-y$ is invertible. Note that the original transformations
$mC_{*},T_{y*}$ and $\tilde{T}_{y*}$  are ring homomorphisms on a point space $\{pt\}$,
and module homomorphisms over any space $X$, by the
multiplicativity with respect to exterior products.
Similarly, these extend by Corollary \ref{cor:mC}.1,
Corollary \ref{cor:Ty}.1 and (\ref{eq:comp}) to transformations 
\begin{equation} \label{eq:compl} \begin{split}
mC^{\wedge}_{*} &: \widehat{M}(var/X)\to G_{0}(X)\otimes \bb{Z}[y][[y^{-1}]]\:,\\
T^{\wedge}_{y*}, \tilde{T}^{\wedge}_{y*} 
&: \widehat{M}(var/X)\to H_{*}(X)\otimes \bb{Q}[y][[y^{-1}]]\:,\\
\end{split} \end{equation}
of the corresponding completions (for $k\to -\infty$)
with respect to the dimension filtration $F_{k}M(var/X)$
of $M(var/X)$. Here
\[F_{k}M(var/X) \quad \text{is generated by} \quad
[X'\to X]\bb{L}^{-n} \quad \text{with $dim(X')-n\leq k$.} \]
\end{rem}

This completion $\widehat{M}(var/X)$ also comes up in the
context of motivic integration. 
In the absolute case $X=\{pt\}$
it was introduced by Kontsevich in his study of the {\em arc-space\/}
${\cal L}(X)$ of $X$ as the value group of a ``motivic measure'' 
$\tilde{\mu}$ on a suitable Boolean algebra of subsets of ${\cal L}(X)$.
This allows one (compare \cite[4.4]{DL} and \cite{V})
 to introduce new invariants for $X$ pure-dimensional, but maybe
singular, as the value of 
\[\tilde{\mu}({\cal L}(X)) \in \widehat{M}(var/\{pt\}) \to R\]
under a suitable homomorphism to a ring $R$.
Instead of $\tilde{\mu}({\cal L}(X))$, one can also use related
``motivic integrals'' over ${\cal L}(X)$.
By our work, one can now introduce similar characteristic classes 
by using a ``relative motivic measure'' $\tilde{\mu}_{X}$ with values in
$\widehat{M}(var/X)$ \cite[sec.4]{Lo}, and the same for ``motivic integrals''
(compare also with \cite{Y6}).  

\begin{ex} \label{ex:motint}
Let $Y$ be a pure-dimensional manifold and $D=\sum_{i=1}^{r}\:a_{i}D_{i}$
be an effective normal crossing divisor (e.g. $a_{i}\in \bb{N}_{0}$) on $Y$, with
smooth irreducible components $D_{i}$. Then one can introduce and evaluate
a motivic integral of the following type (compare \cite{DL,Lo,Craw,V}):
\begin{equation} \label{eq:motint} \begin{split}
\int_{{\cal L}(Y)} \: \bb{L}^{-ord(D)}d\tilde{\mu}_{Y} &=
\sum_{I\subset \{1,\dots,r\}} \: [D_{I}^{o}\to Y] \cdot \prod_{i\in I}\:
\frac{\bb{L}-1}{\bb{L}^{a_{i}+1}-1} \\
&= \sum_{I\subset \{1,\dots,r\}} \: [D_{I}\to Y] \cdot \prod_{i\in I}\:
(\frac{\bb{L}-1}{\bb{L}^{a_{i}+1}-1} -1) \:.
\end{split} \end{equation}
Here we use the notation:
\[D_{I}:=\bigcap _{i\in I}\: D_{i} \quad \text{(with $D_{\emptyset}:=Y$), and}
\quad
D_{I}^{o}:=D_{I}\backslash \ \bigcup_{i\in \{1,\dots,r\}\backslash I} \:D_{i} \:.\]
So by application of a characteristic class $cl_{*}^{\wedge}=
mC^{\wedge}_{*},T^{\wedge}_{y*}, \tilde{T}^{\wedge}_{y*}$ we get:
\begin{equation} \label{eq:motint2} \begin{split}
cl_{*}^{\wedge}\Bigl(\int_{{\cal L}(Y)} \: \bb{L}^{-ord(D)}d\tilde{\mu}_{Y}\Bigr) =
\sum_{I\subset \{1,\dots,r\}} \: cl_{*}([D_{I}^{o}\to Y]) \cdot \prod_{i\in I}\:
\frac{(-y)-1}{(-y)^{a_{i}+1}-1} \\
= \sum_{I\subset \{1,\dots,r\}} \: cl_{*}([D_{I}\to Y]) \cdot \prod_{i\in I}\:
\frac{(-y)-(-y)^{a_{i}+1}}{(-y)^{a_{i}+1}-1}  \:.
\end{split} \end{equation}
\end{ex}
But $D_{I}$ is a closed smooth submanifold of $Y$ so that $cl_{*}([D_{I}\to Y])$ is just the 
pushforward to $Y$ of the corresponding characteristic (homology) class
\[cl_{*}(D_{I})=cl^{*}(TD_{I})\cap [D_{I}] 
\quad \text{for} \quad cl_{*}= mC_{*},T_{y*}, \tilde{T}_{y*}\:.\]

Here the factor $(\bb{L}^{a_{i}+1}-1)^{-1}=\bb{L}^{-(a_{i}+1)}\cdot (1-\bb{L}^{-(a_{i}+1)})^{-1}$ in (\ref{eq:motint}) has to be developed
as the corresponding geometric series in $\widehat{M}(var/\{pt\})$, and similarly for
the factor
$((-y)^{a_{i}+1}-1)^{-1} \in \bb{Z}[y][[y^{-1}]]$ in (\ref{eq:motint2}).
Moreover one gets the last equality in (\ref{eq:motint}) by multiplying out
the following products:
\begin{equation} \label{eq:motint3} \begin{split}
& \prod_{i=1}^{r}\:\Bigl( b_{i}\cdot[D_{i}\to Y]+ [Y\backslash D_{i}\to Y]\Bigr) =\\
& \prod_{i=1}^{r}\: \Bigl( (b_{i}-1)\cdot [D_{i}\to Y]+ [id_{Y}] \Bigr) \in \widehat{M}(var/Y)\:,
\end{split} \end{equation}
with $b_{i}:=(\bb{L}-1)(\bb{L}^{a_{i}+1}-1)^{-1} \in \widehat{M}(var/\{pt\})$.
Recall that multiplication in $\widehat{M}(var/Y)$ is induced from taking the
fiber product over $Y$.\\

Let now $X$ be a normal pure-dimensional algebraic variety which is {\em Gorenstein\/}, i.e., such that the canonical divisor $K_{X}$ is a Cartier divisor. Take a resolution
of singularities $\pi: Y\to X$ such that the ``discrepancy divisor'' 
$D:=K_{Y}- \pi^{*}K_{X}$ is a divisor with normal crossing. Assume $X$ has only {\em log-terminal singularities\/}, i.e. this $D$ is for one (and then for any) such
resolution an effective divisor (so that $X$ has already {\em canonical singularities\/},
since for simplicity we do not consider $\bb{Q}$-divisors). Then one of the main results from motivic
integration says (as an application of the ``transformation rule''), that
the following ``stringy invariant'' of $X$: 
\begin{equation} \label{eq:logter}
{\cal E}_{st}(X)=
\pi_{*}\Bigl(\int_{{\cal L}(Y)} \: \bb{L}^{-ord(D)}d\tilde{\mu}_{Y}\Bigr) \in \widehat{M}(var/X)
\end{equation}
does not depend on the choice of the resolution (compare with \cite[(7.7),(7.8)]{V}
and \cite{DL,Lo,Craw}).
So this is an intrinsic invariant of $X$, and we can introduce the corresponding
{\em stringy characteristic homology class\/} $cl_{*}^{st}(X)$ of $X$ for $cl_{*}=
mC_{*},T_{y*}, \tilde{T}_{y*}$ by
\begin{equation} \label{eq:stringcl} \begin{split}
cl_{*}^{st}(X):&=
cl_{*}^{\wedge}\pi_{*}\Bigl(\int_{{\cal L}(Y)} \: \bb{L}^{-ord(D)}d\tilde{\mu}_{Y}\Bigr) \\
&= \pi_{*}cl_{*}^{\wedge}\Bigl(\int_{{\cal L}(Y)} \: \bb{L}^{-ord(D)}d\tilde{\mu}_{Y}\Bigr) \:.
\end{split} \end{equation}

Then the {\em stringy Hirzebruch classes\/} $T^{st}_{y*}(X)$ and $\tilde{T}^{st}_{y*}(X)$
interpolate in the following sense between the {\em elliptic class\/}
${\cal E}ll(X)$ of $X$ defined
by Borisov-Libgober \cite{BL1,BL2} and the {\em stringy $E$-function\/} $E_{st}(X)$
of Batyrev \cite{Bat} for $X$ complex algebraic and complete:
\begin{equation} \label{eq:ell=hirze}
\lim_{\tau\to i\infty} \: {\cal E}ll(X)(z,\tau) = y^{-1/2\cdot dim\;X}\cdot
\tilde{T}^{st}_{-y*}(X) \quad \text{for $y=e^{2\pi iz}$} \:,
\end{equation}
and
\begin{equation} \label{eq:Efct} \begin{split}
\chi_{-y}^{st}(X):&= \int_{[X]}\: T^{st}_{-y*}(X)\\
&= \int_{[X]}\: \tilde{T}^{st}_{-y*}(X)  = E_{st}(X)(y,1)
\end{split} \end{equation}
for $X$ (complex algebraic and) complete. Here the equality (\ref{eq:ell=hirze})
for the {\em elliptic class\/} of $X$
in the sense of Borisov-Libgober \cite[def.3.2,rem.3.3]{BL2},
${\cal E}ll(X):={\cal E}ll_{orb}(X,E,G)$  with $G=\{id\}$ and $E=\emptyset$,
follows directly from (\ref{eq:motint2}) and the calculation given in
\cite[p.14/15]{BL2}. Similarly, (\ref{eq:stringcl}) follows from (\ref{eq:motint2})
and the definition of $E_{st}(X)(u,v):=E_{st}(X,\emptyset)(u,v)$ given in \cite[def.3.7]{Bat}, if one takes
$(u,v)=(y,1)$ (compare with (\ref{eq:E}) in Section \ref{sec:Hodge}).\\

So these stringy Hirzebruch classes are ``in between'' the elliptic class
and the stringy $E$-function, and as suitable limits they are ``weaker''
than these more general invariants. But they have the following good properties of both of them:
\begin{itemize}
\item The stringy Hirzebruch classes come from a functorial ``additive''
characteristic homology class.
\item The stringy $E$-function comes from the ``additive'' {\em $E$-polynomial\/} 
defined by Hodge theory on
$K_{0}(var/\{pt\})$, which doesn't have a homology class version
(compare with section \ref{sec:Hodge}).
\item The elliptic class is a homology class, which doesn't come
from an ``additive'' characteristic class (of vector bundles), since the corresponding
{\em elliptic genus\/} is more general than the {\em Hirzebruch $\chi_{y}$-genus\/},
which is the most general ``additive'' genus of such a class
(as explained in the introduction).
\end{itemize}
 
Finally the stringy Hirzebruch class $T^{st}_{y*}(X)$ specializes for $y=-1$
in the following way to the {\em stringy Chern class\/} $c_{*}^{st}(X)$ of $X$ as introduced in \cite{A2,FLNU}:
\begin{equation} \label{eq:hirze=ch}
\lim_{y\to -1}\:T^{st}_{y*}(X) = c_{*}^{st}(X) \in A_{*}(X)\otimes\bb{Q} \:.
\end{equation}
This follows directly from (\ref{eq:motint2}) and the left square of the commutative
diagram (\ref{eq:real1}), if one compares this with 
\cite[sec.5.5,6.5]{A2} and \cite[sec.3,def.4.1]{FLNU}.\\

If we specialize in (\ref{eq:motint2}) to $y=0$, then we get by ``additivity'':
\begin{equation} \label{eq:crep1}
\lim_{y\to 0}\: T^{st}_{y*}(X) = \pi_{*}(Td^{*}(TY)\cap [Y])
\end{equation}
so that the welldefinedness of this limit is just a special case of
Corollary \ref{cor:tdbirat}. Assume now that $\pi: Y\to X$ is a {\em crepant resolution\/},
i.e. all multiplicities $a_{i}$ of the ``discrepancy divisor'' $D$ are zero.
Then one gets in the same way:
\begin{equation} \label{eq:crep2}
T^{st}_{y*}(X) = \pi_{*}(T_{y}^{*}(TY)\cap [Y])
\end{equation}
and
\begin{equation} \label{eq:crep3}
c_{*}^{st}(X) =
\lim_{y\to -1}\: T^{st}_{y*}(X) = \pi_{*}(c^{*}(TY)\cap [Y])\:.
\end{equation}
In particular the right hand side does not depend on the choice of a
crepant resolution (compare \cite[cor.1.2]{A1} and \cite[prop.4.4]{FLNU}).

\section{Comparison with functorial $L$-classes}  
\label{sec:Lclass}

In this section we work in the algebraic context over the base field $k=\bb{C}$,
or in the complex analytic context with compact spaces,
and explain the relation of our {\em motivic class\/} $T_{1*}$ to
the {\em $L$-class transformation\/} of Cappell-Shaneson \cite{CS1}.\\

Let $X$ be a complex analytic (algebraic) space with $A:=D^{b}_{c}(X)$ the bounded derived
category of complex analytically (algebraically) constructible complexes of sheaves of
$\bb{K}$-vector spaces (for $\bb{K}$ a subfield of $\bb{R}$, compare \cite{KS,Sch1}).
So we consider bounded sheaf complexes ${\cal F}$, which have locally constant cohomology
sheaves with finite dimensional stalks along the strata of a complex analytic (algebraic)
Whitney stratification $X_{\bullet}$ of $X$. In the algebraic context, or for a compact analytic
space $X$, such a stratification has only finitely many strata.
Then $A$ is a triangulated category with translation
functor $T_{A}=T=[1]$ and a duality $D_{A}$ in the sense of \cite{You,Ba1,Ba2,GN} induced by
the {\em Verdier duality functor\/} (compare \cite[chap.4]{Sch1} and \cite[chap.VIII]{KS}):
\[D_{X}:= Rhom(\cdot, k^{!}\bb{K}_{pt}): \: D^{b}_{c}(X)\to D^{b}_{c}(X)\:,\]
with $k:X\to \{pt\}$ a constant map. Here a {\em triangulated category $B$ with translation
functor $T_{B}$ has a duality $D_{B}$\/} in the sense of \cite{You,Ba1,Ba2,GN}, if
\begin{enumerate}
\item $D_{B}: B\to B$ is a contravariant functor with $D_{B}\circ T_{B}=T_{B}^{-1}\circ D_{B}$,
which preserves distinguished triangles,
\item one has a {\em biduality isomorphism\/} $can: id\stackrel{\sim}{\to} D_{B}\circ D_{B}$ such that $can_{T_{B}(M)}=T_{B}(can_{M})$ (i.e. $can$ commutes with translation)
and $D_{B}(can_{M})\circ can_{D_{B}(M)}$ $= id_{D_{B}(M)}$
for any $M\in ob(B)$.
\end{enumerate}
For the biduality isomorphism $can$ in the case of the 
Verdier duality functor $D_{X}$ compare with \cite[cor.4.2.2]{Sch1} and \cite[prop.3.4.3,prop.8.4.9]{KS}.

\begin{rem} \label{rem:dual}
Our notion of a duality $D_{B}$ in a triangulated category $B$ corresponds to a
{\em $\delta=1$-duality\/} in the sense of \cite{Ba1,Ba2,GN}. These authors also consider a {\em $\delta=-1$-duality\/} corresponding to the case that $D_{B}$ maps distinguished triangles to the negative of distinguished triangles. 
\end{rem}

A constructible complex ${\cal F}\in ob(D^{b}_{c}(X))$ is called {\em selfdual\/}
(and similarly in the more general context of a triangulated category with duality),
if there is an isomorphism
\[ d: {\cal F} \stackrel{\sim}{\to} D_{X}({\cal F}) \:.\]
The pair $({\cal F},d)$ is called {\em symmetric\/} or {\em skew-symmetric\/} 
(in \cite{You} this is called
``even'' or ``odd''), if
\[ D_{X}(d)\circ can = d \quad \text{or} \quad D_{X}(d)\circ can = -d \:.\]
Note that a skew-symmetric pair $({\cal F},d)$ for the biduality isomorphism $can$
is just the same as a symmetric pair $({\cal F},d)$ for the biduality isomorphism $-can$.\\

Finally an isomorphism or {\em isometry\/} of selfdual objects $({\cal F},d)$ and $({\cal F}',d')$
is an isomorphism $u$ such that the following diagram commutes:
\begin{displaymath} \begin{CD}
{\cal F} @> u > \sim > {\cal F}' \\
@V d VV @VV d' V \\
D_{X}({\cal F}) @< \sim < D_{X}(u) < D_{X}({\cal F}') \:.
\end{CD} \end{displaymath}

We consider now the complex algebraic context, or a compact analytic space $X$
so that the isomorphism classes of such (skew-)symmetric selfdual complexes form a set.
This becomes a {\em monoid\/} with addition induced by the direct sum.
Then the {\em Witt group\/} $W_{+}(X)$ ($W_{-}(X)$, resp.) of such symmetric 
(skew-symmetric, resp.) selfdual complexes on $X$ is the
quotient of this monoid with respect to the submonoid of {\em neutral\/}
symmetric (skew-symmetric, resp.) selfdual complexes in the
sense of \cite{Ba1,Ba2}. Here $({\cal F},d)$ is called {\em neutral\/}, if there exists an isomorphism
of distinguished triangles (for some ${\cal L} \in ob(D^{b}_{c}(X))$):

\begin{equation} \label{eq:neutral} \begin{CD}
D_{X}({\cal L})[-1] @> w >> {\cal L} @> \alpha >> {\cal F} @> \beta >> D_{X}({\cal L}) \\
@|  @V can V \wr V @VV d V @| \\
D_{X}({\cal L})[-1] @>> D_{X}(w)[-1] > D_{X}D_{X}({\cal L}) @>> D_{X}(\beta) > D_{X}({\cal F}) 
@>>  D_{X}(\alpha) > D_{X}({\cal L}) \:.
\end{CD} \end{equation}

Then $W_{\pm}(X)$ is indeed an abelian group, since $({\cal F},d)\oplus ({\cal F},-d)$
is always neutral! These Witt groups in the sense of Balmer are different (!) from a corresponding notion 
introduced by Youssin \cite{You}, based on his notion of an {\em elementary cobordism\/} in the context of
a triangulated category with duality. We will call these groups {\em cobordism groups\/} $\Omega_{\pm}(X)$
of selfdual constructible complexes on $X$. Youssin's starting point is the notion of an {\em octahedral diagram\/}
in a triangulated category with duality $D$, i.e.  a diagram $(Oct)$ of the following form:

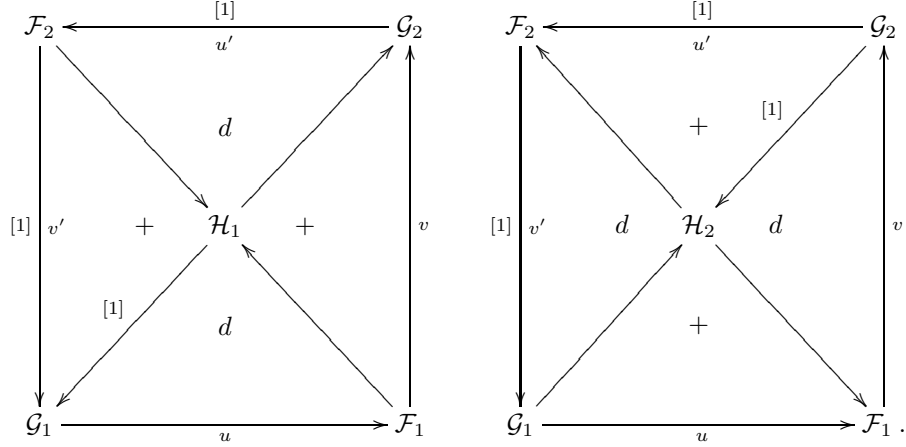
\xymatrix{{\cal F}_{2} \ar[dddd]_{[1]}^{v'} \ar[ddr] & & 
{\cal G}_{2} \ar[ll]^{u'}_{[1]} &  {\cal F}_{2} \ar[dddd]_{[1]}^{v'}
 & &  {\cal G}_{2} \ar[ll]^{u'}_{[1]} \ar[ddl]_{[1]} \\
& d & & & + \\
& + \quad\quad {\cal H}_{1} \quad\quad + \ar[ddl]_{[1]} \ar[uur] & & &
 d\quad\quad {\cal H}_{2}\quad\quad d \ar[uul] \ar[ddr] \\
& d & & & + \\
{\cal G}_{1} \ar[rr]_{u} & & {\cal F}_{1} \ar[luu] \ar[uuuu]_{v} 
& {\cal G}_{1} \ar[rr]_{u} \ar[ruu]
 & & {\cal F}_{1}  \ar[uuuu]_{v} \:.}

Here the morphism marked by $[1]$ are of degree one, the triangles marked $+$ are commutative,
and the ones marked $d$ are distinguished. Finally the two composite morphisms from ${\cal H}_{1}$
to ${\cal H}_{2}$ (via ${\cal G}_{1}$ and ${\cal G}_{2}$) have to be the same, and similarly for
the two composite morphisms from ${\cal H}_{2}$
to ${\cal H}_{1}$ (via ${\cal F}_{1}$ and ${\cal F}_{2}$).\\

Application of the duality functor $D$ and a rotation by $180^{o}$ 
about the axis connecting upper-left and lower-right corner induces
another octahedral diagram $(RD\cdot Oct)$ such that $RD$ applied to $(RD\cdot Oct)$
gives the octahedral diagram $(D^{2}\cdot Oct)$ which one gets from $(Oct)$ by application of
$D^{2}$ (compare with \cite[p.387/388]{You} for more details).
Then the octahedral diagram $(Oct)$ is called {\em symmetric\/}
or {\em skew-symmetric\/} (in \cite{You} this is called  ``even'' or ``odd''), if there is an isomorphism $d: (Oct)\to (RD\cdot Oct)$ of
octahedral diagrams such that
\[ RD(d)\circ can = d \quad \text{or} \quad RD(d)\circ can = -d \]
as maps of octahedral diagrams $(Oct)\to (RD\cdot Oct)$
(compare \cite[def.6.1]{You}). Note that this induces in particular 
(skew-)symmetric dualities $d_{1}$ and $d_{2}$ of the corners ${\cal F}_{1}$ and ${\cal F}_{2}$,
and $(Oct,d)$ is called an {\em elementary cobordism\/} between $({\cal F}_{1},d_{1})$ and
$({\cal F}_{2},d_{2})$. This notion is a symmetric relation \cite[rem.6.2]{You}.
Similarly, $({\cal F},d)$ is elementary cobordant to itself
(use the octahedral diagram $(Oct)$ with 
${\cal F}_{1}={\cal F}_{2}={\cal H}_{1}={\cal H}_{2}$, ${\cal G}_{1}=0={\cal G}_{2}$
and the corresponding isomorphism induced by $d$).
$({\cal F},d)$ and $({\cal F}',d')$ are called {\em cobordant\/} (compare \cite[p.528]{CS1}),
if there is a sequence 
\[({\cal F},d)=({\cal F}_{0},d_{0}),\: ({\cal F}_{1},d_{1}),\:\dots \:, 
({\cal F}_{m},d_{m})=({\cal F}',d')\]
with $({\cal F}_{i},d_{i})$ elementary cobordant to $({\cal F}_{i+1},d_{i+1})$ for 
$i=0,\dots,m-1$. This {\em cobordism relation\/} is then an equivalence relation.\\

The {\em cobordism group\/} $\Omega_{+}(X)$ ($\Omega_{-}(X)$, resp.)
of selfdual constructible complexes on $X$ is the quotient of the monoid of isomorphism classes
of symmetric (skew-symmetric, resp.) selfdual complexes by this
cobordism relation. These are again monoids, since this relation commutes with direct sums.
But a neutral selfdual constructible complex $({\cal F},d)$ in the sense of Balmer given by a 
distinguished triangle (\ref{eq:neutral}) induces a (skew-)symmetric octahedral diagram
of the following type (with isomorphism $d$ induced by the isomorphism of this distinguished triangle):

\xymatrix{0 \ar[dddd]^{[1]} \ar[ddr] & & 
D({\cal L}) \ar[ll]_{[1]} &  0 \ar[dddd]^{[1]}
 & &  D({\cal L}) \ar[ll]_{[1]} \ar[ddl]_{[1]} \\
& d & & & + \\
& + \quad\quad D({\cal L}) \quad\quad + \ar[ddl]_{[1]} \ar@{=}[uur] & & &
 d\quad\quad {\cal L}\quad\quad d \ar[uul] \ar[ddr]_{\alpha} \\
& d & & & + \\
{\cal L} \ar[rr]_{\alpha} & & {\cal F} \ar[luu]_{\beta} \ar[uuuu]^{\beta} 
& {\cal L} \ar[rr]_{\alpha} \ar@{=}[ruu]
 & & {\cal F}  \ar[uuuu]^{\beta} \:.}

So a neutral selfdual constructible complex $({\cal F},d)$ is elementary corbordant to $0$.
But then $\Omega_{\pm}(X)$ is also an abelian group together with a canonical group epimorphism
(compare \cite[introduction]{Ba1} and \cite[rem.3.25]{Ba2}):
\[W_{\pm}(X) \to \Omega_{\pm}(X) \to 0\:.\]

Consider now a proper algebraic (or holomorphic) map $f: X\to Y$,
with $X,Y$ compact in the analytic context. 
Then $Rf_{*}\simeq Rf_{!}$ maps $D^{b}_{c}(X)$ to $D^{b}_{c}(X)$
(compare \cite[chap.4]{Sch1} and \cite[chap.VIII]{KS}).
Moreover, the {\em adjunction isomorphism\/} (\cite[p.120]{Sch1}, \cite[prop.3.1.10]{KS}):
\[ Rf_{*}Rhom({\cal F},f^{!}k^{!}\bb{K}_{pt}) \simeq
Rhom(Rf_{!}{\cal F},k^{!}\bb{K}_{pt}) \]
induces the isomorphism
\begin{equation} \label{eq:dual1}
Rf_{*}D_{X} \stackrel{\sim}{\to} D_{Y}Rf_{!} \simeq  D_{Y}Rf_{*} 
\end{equation}
so that $Rf_{*}$ {\em commutes with Verdier-duality\/}
(compare with \cite[def.1.8]{GN} for the abstract notion of a {\em duality preserving functor\/}
between triangulated categories with duality). In particular $Rf_{*}$ maps 
selfdual constructible complexes on $X$ to selfdual constructible complexes
on $Y$ inducing group homomorphisms
\[f_{*}: \Omega_{\pm}(X)\to \Omega_{\pm}(Y);\: [({\cal F},d)]\mapsto [(Rf_{*}{\cal F},Rf_{*}(d))] \:.\]

Before we can compare this with the corresponding notion from \cite{CS1,Y1},
we have to explain the relation between {\em duality\/} $D=D_{X}$ and {\em translation\/} $T=[1]$. 
By the equality $D_{X}\circ T = T^{-1}\circ D_{X}$ one gets the {\em translated duality\/}
$T^{2n}D_{X}:= T^{2n}\circ D_{X}$ ($n\in \bb{Z}$) with the biduality isomorphism
\[ can: id \stackrel{\sim}{\to} D_{X}\circ D_{X} \simeq D_{X}\circ T^{-2n}\circ T^{2n}\circ D_{X}
\simeq T^{2n}D_{X}\circ T^{2n}D_{X} \:.\]
Then the translation $T^{n}$ ($n\in \bb{Z}$) induces an isomorphism of {\em cobordism groups\/}
(and similarly for {\em Witt groups\/})
\begin{equation} \label{eq:shift} \begin{split}
T^{n}: \Omega(X):=\Omega(X,D_{X}) & :=\Omega_{+}(X,D_{X})\oplus \Omega_{-}(X,D_{X}) 
\stackrel{\sim}{\to} \Omega(X,T^{2n}D_{X}); \\
& [({\cal F},d)]\mapsto [(T^{n}({\cal F}),T^{n}(d))] \:, 
\end{split} \end{equation}
with
\[T^{n}(d): T^{n}({\cal F}) \stackrel{\sim}{\to} T^{n}\circ D_{X}({\cal F})
= T^{2n}\circ T^{-n}\circ D_{X}({\cal F}) = T^{2n}D_{X}(T^{n}({\cal F})) \:.\]

Note that this isomorphism is {\em parity\/} preserving (or changing) depending on $n$ even (or odd).
This follows from the {\em anti-commutative\/} diagram (compare \cite[rem.1.10.16]{KS}):
\begin{displaymath} \begin{CD}
Hom({\cal F}[1],{\cal G}[1]) @>\sim >>  Hom({\cal F},{\cal G}[1])[-1] \\
@V \wr VV  @VV \wr V \\
Hom({\cal F}[1],{\cal G})[1]  @>>\sim > Hom({\cal F},{\cal G}) \:.
\end{CD} \end{displaymath}

\begin{rem} We only need to consider even shifts $T^{2n}D_{X}$ of the Verdier duality.
Note that $T$ maps by definition distinguished triangles to the negative of distinguished
triangles so that $T\circ D_{X}$ is a {\em $\delta=-1$-duality\/} in the sense of \cite{Ba1,Ba2,GN}.
\end{rem}

In the following we identify the shifted cobordism groups by the isomorphism (\ref{eq:shift}).
Note that \cite{CS1,Y1} use the {\em shifted duality\/} $T^{2n(X)}D_{X}$ with $n(X)$ the complex dimension of $X$
in their definition of selfdual constructible complexes. Then the following diagram becomes
commutative under these identifications:
\begin{displaymath} \begin{CD}
\Omega(X) @> T^{n(X)} > \sim > \Omega(X,T^{2n(X)}D_{X}) \\
@V f_{*} VV  @VV f_{*} V \\
\Omega(Y) @> T^{n(Y)} > \sim > \Omega(Y,T^{2n(Y)}D_{Y}) \:,
\end{CD} \end{displaymath}
with $f_{*}$ on the right hand side induced by $Rf_{*}[n(Y)-n(X)]$ as in \cite[prop.4.1]{CS1} and
\cite[prop.1.6]{Y1}.\\

Let us now explain an important example of a (skew-)symmetric selfdual
constructible complex.

\begin{ex} \label{ex:cob}
Assume $Z$ is a complex manifold of pure dimension $n$ so that the {\em complex orientation\/} of $Z$ induces an isomorphism
$k^{!}\bb{K}_{pt} \simeq \bb{K}_{Z}[2n]$ (compare \cite[sec.III.3.3]{KS}).
Let ${\cal L}$ be a {\em Poincar\'e local system\/} on $Z$, i.e. a locally constant sheaf of finite rank
with a symmetric nondegenerate bilinear pairing
\[\phi: {\cal L}\times {\cal L} \to \bb{K}_{Z}\:.\]
For example we can take ${\cal L}=\bb{K}_{Z}$ with the obvious pairing given by multiplication.
Then $\phi$ induces an isomorphism
\[d': {\cal L} \stackrel{\sim}{\to} hom({\cal L},\bb{K}_{Z}) \simeq
Rhom({\cal L},\bb{K}_{Z})= D_{Z}({\cal L})[-2n] \:.\]
In particular $d=d'[n]$ induces a (skew-)symmetric selfduality of ${\cal L}[n]$.

Assume moreover that $Z$ is a locally closed constructible subset of the analytic space $X$,
with $i: \bar{Z}\to X$ the closed inclusion of the complex analytic closure $\bar{Z}$ of $Z$.
Then $d$ induces a (skew-)symmetric selfduality of the {\em twisted intersection cohomology complex\/}
\[ {\cal IC}^{\bullet}_{\bar{m}}(\bar{Z},{\cal L})[n] \] 
on $\bar{Z}$ (compare \cite[sec.6.02]{Sch1}, with $\bar{m}$ the middle perversity).
Here we use the convention that ${\cal IC}^{\bullet}_{\bar{m}}(\bar{Z},{\cal L})$
gives back ${\cal L}$ by restriction to the open set $Z$ of $\bar{Z}$.
By proper pushdown along the map $i$ we finally get the (skew-)symmetric selfdual
complex $Ri_{*}{\cal IC}^{\bullet}_{\bar{m}}(\bar{Z},{\cal L})[n]$ on $X$.
\end{ex}

We also have to point out that \cite{CS1,Y1} use a slightly {\em weaker\/} notion of ``elementary
cobordism'' for the definition of their cobordism groups of selfdual constructible complexes.
In terms of the unshifted Verdier duality functor, this is defined as follows:
Start with morphisms
\begin{displaymath} \begin{CD} 
{\cal G} @> u >> {\cal F} @> v >> {\cal H} 
\end{CD} \end{displaymath} 
in $D^{b}_{c}(X)$ with $v\circ u=0$. Suppose that there is an isomorphism
${\cal H} \simeq D_{X}({\cal G})$ such that the following diagram commutes:
\begin{displaymath} \begin{CD}
{\cal F} @> v >> {\cal H} \\
@V \wr V d V @VV \wr V \\
D_{X}({\cal F}) @>> D_{X}(u) > D_{X}({\cal G}) \:.
\end{CD} \end{displaymath}
Then the morphisms $u$ and $v$ can be included into a octahedral diagram $(Oct)$
with $X_{2}=:C_{u,v}$ an ``iterated cone'' in the sense of \cite{CS1}
(compare \cite[lem.5.3]{You}). Moreover one can lift the duality $d: {\cal F}
\stackrel{\sim}{\to} D_{X}({\cal F})$ to a duality
\[d: C_{u,v} \stackrel{\sim}{\to} D_{X}(C_{u,v}) \:,\]
and $({\cal F},d)$ is by their definition elementary cobordant to $(C_{u,v},d)$.\\

But the induced duality of $C_{u,v}$ is {\em non-canonical\/} and therefore also {\em non-functorial\/},
since ``the cone'' of a morphism is not canonically defined in a triangulated category.
In particular, it is not clear in general, if the corresponding octahedral diagram
$(Oct)$ can be chosen to be selfdual and (skew-)symmetric (compare \cite[p.390]{You})!
This is indeed the case if $Hom({\cal G}[1],{\cal H})=0$, e.g. if ${\cal G}$ and
${\cal H}$ are perverse sheaves with respect to the middle perversity $t$-structure (compare 
\cite[chapter X]{KS}). On the other hand, it is clear by definition, that an elementary
cobordism in the sense of \cite{You} is functorial, and that it induces an elementary
cobordism in the sense of \cite{CS1}. Therefore the results of \cite[sec.5]{CS1}
can be reformulated as in \cite[cor.2.3]{Y1}:

\begin{thm}[Cappell-Shaneson] \label{thm:Lclass}
For a compact complex analytic (or algebraic) space $X$ there is a {\em homology $L$-class transformation\/}
\[L_{*}: \Omega(X) =  \Omega_{+}(X)\oplus \Omega_{-}(X) \to H_{2*}(X,\bb{Q}) \]
as in (\ref{eq:real2}), which is a group homomorphism functorial for the pushdown $f_{*}$
induced by a holomorphic (or algebraic) map. 
The degree of $L_{0}(({\cal F},d))$ is the {\em signature\/} of the induced pairing 
\[H^{0}(X,{\cal F})\otimes_{K}\bb{R} \times H^{0}(X,{\cal F})\otimes_{K}\bb{R} \to \bb{R} \]
(by definition this is $0$ for a skew-symmetric pairing).
Moreover, for $X$ smooth of pure dimension $n$ one has the normalization 
\[L_{*}((\bb{K}_{X}[n],d)) = L^{*}(TX)\cap [X] \:,\]
with $(\bb{K}_{X}[n],d)$ as in example \ref{ex:cob}.
\hfill $\Box$
\end{thm}

If $h:X'\to X$ is an isomorphism of purely n-dimensional complex manifolds, then
\[ h_{*}((\bb{K}_{X'}[n],d)) = (\bb{K}_{X}[n],d)\]
(with $d$ as in example \ref{ex:cob}). Consider now the complex algebraic context,
or the analytic context with $X$ compact. 
Then there is a unique group homomorphism 
\begin{equation} \label{eq:defsd}
sd : Iso^{pro}(sm/X) \to \Omega(X)
\end{equation}
satisfying the {\em normalization\/} condition 
\[sd([f:X'\to X])=f_{*}\bigl( [(\bb{K}_{X'}[n],d)] \bigr)  \]
for $X'$ purely n-dimensional. But contrarily to what is claimed in \cite{Y1},
this transformation $sd$ is in general not surjective!\\

Before we explain this, we have to recall another important result from
\cite{CS1,You}. Let
\[(^{m}D_{c}^{\leq 0}(X),\,^{m}D_{c}^{\geq 0}(X)) \]
be the {\em perverse t-structure\/} on $D_{c}^{b}(X)$ with respect to the middle perversity $m$.
(compare \cite[chapter X]{KS} and \cite[chapter 6]{Sch1}). The Verdier duality $D_{X}$ interchanges $^{m}D_{c}^{\leq 0}(X)$ and $^{m}D_{c}^{\geq 0}(X)$
\cite[prop.10.3.5]{KS} so that the corresponding {\em perverse cohomology functor\/} 
$^{m}{\cal H}^{0}=\tau^{\geq 0}\tau^{\leq 0}$ \cite[def.10.1.9]{KS} commutes with Verdier duality. 
Then one has by \cite[Ex.6.6]{You} (and compare with \cite[lem.3.3]{CS1}) the important

\begin{lem}[Youssin] \label{lem:perverse}
Assume $({\cal F},d)$ is a (skew-)symmetric selfdual
constructible complex on $X$. Then $({\cal F},d)$ is elementary cobordant
to the (skew-)\-symmetric selfdual constructible complex
$(^{m}{\cal H}^{0}({\cal F}),\,^{m}{\cal H}^{0}(d))$.
\hfill $\Box$
\end{lem}

Consider a proper holomorphic map 
$f:X'\to X$ of pure-dimensional complex manifolds such that 
all higher direct image sheaves $R^{i}f_{*}\bb{K}_{X'}$ ($i\in \bb{Z}$) are {\em locally constant\/}.
Let $n(X),n(X')$ be the complex dimensions of $X$ and $X'$. Then
\begin{equation} \label{eq:pervdirect}
^{m}{\cal H}^{0}(Rf_{*}\bb{K}_{X'}[n(X')]) \simeq R^{i}f_{*}\bb{K}_{X'}[n(X)]
\quad \text{for $i=n(X)-n(X')$.}
\end{equation}
Note that $i$ is just the complex fiber dimension of $f$.
Of course, in general the higher direct image sheaves $R^{i}f_{*}\bb{K}_{X'}$ are only {\em constructible\/}
for such a proper holomorphic map $f$.
But assume that $X=\bb{P}^{1}(\bb{C})$ so that 
\[f: f^{-1}(D^{*}_{r}(0))\to D^{*}_{r}(0):=\{z\in \bb{C}|\: 0<|z|<r\} \] 
satisfies the assumption above for $r$ small enough. Then the {\em monodromy automorphism\/} 
of the complexification $(\cdot)\otimes_{K} \bb{C}$ of the local system
\[ ^{m}{\cal H}^{0}(Rf_{*}\bb{K}_{X'}[n(X')])|D^{*}_{r}(0) \]
has only {\em roots of unity as eigenvalues\/} by the well known ``monodromy theorem''. 
Take now an irreducible Poincar\'e local system ${\cal L}$ on $\bb{C}^{*}$, whose corresponding
monodromy automorphism (of the complexification) has not roots of unity as eigenvalues
(e.g. a suitable rotation of ${\cal L}_{z}:=\bb{C}=\bb{R}^{2}$ for $\bb{K}=\bb{R}$ and
$z\in \bb{C}^{*}$). Then
\[ ({\cal IC}^{\bullet}_{\bar{m}}(\bb{P}^{1}(\bb{C}),{\cal L})[1],d) \not\in 
sd\bigl(Iso^{pro}(sm/\bb{P}^{1}(\bb{C})) \bigr) \:,\]
since this is a {\em simple\/} (selfdual) perverse sheaf, and $\Omega(X)$ is
freely generated by the isomorphism classes of such simple (selfdual) perverse sheaves
\cite[cor.7.5]{You}.\\

The mistake of \cite{Y1} can already be explained for a {\em holomorphic submersion\/}
$f:X'\to X$ of compact pure-dimensional complex manifolds.
Then all higher direct image sheaves $R^{i}f_{*}\bb{K}_{X'}$ ($i\in \bb{Z}$) are locally constant
by the ``Ehresmann fibration theorem'' (compare also with \cite[Ex.4.1.2]{Sch1}). Then
one gets by Lemma \ref{lem:perverse} and (\ref{eq:pervdirect}):
\[f_{*}\bigl([(\bb{K}_{X'}[n(X')],d')]\bigr) = [(R^{i}f_{*}\bb{K}_{X'}[n(X)],d)] \]
for $i=n(X)-n(X')$ and a suitable induced duality $d$.
And in general
\[ [(R^{i}f_{*}\bb{K}_{X'}[n(X')],d)] \neq rk\cdot [(\bb{K}_{X}[n(X)],d)] \]
with $rk$ the {\em rank\/} of the local system $R^{i}f_{*}\bb{K}_{X'}$,
contrarily to what is claimed in \cite[p.1011]{Y1}. 
If $X$ is connected and this local system is already {\em constant\/}, then one gets 
for $\bb{K}=\bb{R}$ by \cite[prop.4.4]{CS1}:
\begin{equation} \label{eq:sign1}
[(R^{i}f_{*}\bb{R}_{X'}[n(X')],d)] = sign(F_{y})\cdot [(\bb{R}_{X}[n(X)],d)] \:,
\end{equation}
with $sign(F_{y})$ the {\em signature of the fiber\/} $F_{y}:=f^{-1}(\{y\})$ for $y\in X$.
In particular theorem \ref{thm:Lclass} implies in this case the ``multiplicativity theorem'' of Chern, Hirzebruch and Serre (compare \cite[p.42]{HBJ}):
\[sign(X')= sign(F_{y})\cdot sign(X) \:.\]

But if the corresponding local system $R^{i}f_{*}\bb{K}_{X'}$ is {\em not\/} constant,
then one gets in general no equality of the form
\[ [(R^{i}f_{*}\bb{K}_{X'}[n(X')],d)] = m\cdot [(\bb{K}_{X}[n(X)],d)] \]
for some $m\in \bb{Z}$. Consider for example a compact complex (algebraic)
surface $X'$ with $sign(X')\neq 0$,
which fibers by a holomorphic (algebraic) submersion over a compact connected holomorphic curve $X$.
Then
\[sign(X')=deg\: L_{0}\bigl([(R^{1}f_{*}\bb{K}_{X'}[2],d)]\bigr)
\neq m\cdot deg\: L_{0}\bigl([(\bb{K}_{X}[1],d)]\bigr) =0 \:.\]
Examples of such surfaces are due to Atiyah \cite{At} and Kodaira.\\

Now we are ready for the main result of this section.

\begin{thm} \label{thm:sd}
Consider the complex algebraic context. Then the {\em selfduality transformation\/} 
$sd : Iso^{pro}(sm/X) \to \Omega(X)$ satisfies the ``blow-up relation''
(\ref{eq:bl}) of Theorem \ref{thm:blowup} for smooth pure-dimensional spaces with
$f=id_{X}$. It induces therefore by Corollary \ref{cor:bl} a unique group homomorphism
\[sd : K_{0}(var/X) \to \Omega(X)\]
commuting with proper pushdown and satisfying the {\em normalization condition\/}
\[sd([id_{X}])=[(\bb{K}_{X}[n],d)]  \]
for $X$ smooth and purely n-dimensional. The same is true in
the analytic context for 
\[sd : K_{0}(an/X) \to \Omega(X)\]
with $X$ compact.
\end{thm}

\begin{proof}
Let us consider the blow-up diagram
\begin{displaymath} \begin{CD}
E @> i'>> Bl_{Y}X=X' \\
@VV q' V  @VV q V \\
Y @> i >> X  \:,
\end{CD} \end{displaymath}
with $i$ a closed embedding of smooth pure dimensional spaces.
Here $Bl_{Y}X\to X$ is the blow-up of $X$ along 
$Y$ with exceptional divisor $E$. 
Let $d(E),d(X),d(Y)$ be the complex dimension of the corresponding manifolds,
with $m:=d(X)-d(Y)$ the complex codimension of $Y$ in $X$. 
The case $m=0$ is obvious, since $Y$ is a union of irreducible components of $X$, with $E=\emptyset$ and
$X'=X\backslash Y$. So we can assume $m>0$.
Then
\[q': E=\bb{P}(N_{Y}X)\to Y\]
is just the projection of the projective bundle $\bb{P}(N_{Y}X)$ corresponding
to the normal bundle $N_{Y}X$ of $Y$ in $X$. 
In particular all higher direct image sheaves
$R^{k}q'_{*}(\bb{K}_{E})$ are constant:
\[ R^{k}q'_{*}(\bb{K}_{E}) =
\begin{cases}
\bb{K}_{Y} &\text{for $k=0,2,\dots,2m-2$,}\\
0 &\text{otherwise.}
\end{cases} \]
Moreover this is already true for $\bb{K}=\bb{Z}$ so that
\[ q'_{*}\bigl([(\bb{K}_{E}[n(E)],d)]\bigr) = \epsilon \cdot  [(\bb{K}_{Y}[n(Y)],d)] \]
for a suitable (locally constant) $\epsilon \in \{-1,0,1\}$. But then we can assume $\bb{K}=\bb{R}$
so that by Lemma \ref{lem:perverse} and (\ref{eq:sign1}):
\begin{equation} \label{eq:proof1}
q'_{*}\bigl([(\bb{K}_{E}[n(E)],d)]\bigr) = sign(\bb{P}^{m-1}(\bb{C}))\cdot  [(\bb{K}_{Y}[n(Y)],d)] \:.
\end{equation}

Let $j:X\backslash Y \to X$ be the open inclusion of the complement of $Y$, and consider
the distinguished triangle:
\begin{displaymath} \begin{CD}
\bb{K}_{X} @> ad_{q} >> Rq_{*}\bb{K}_{X'} @>>> {\cal K} @> [1] >> \:.
\end{CD} \end{displaymath}

Then $ad_{q}: \bb{K}_{X} \stackrel{\sim}{\to} R^{0}q_{*}\bb{K}_{X'}$ is an isomorphism,
with $j^{*}{\cal K} =0$ so that $Ri_{*}i^{*}{\cal K}\simeq {\cal K}$.
Pulling back this triangle by $i^{*}$ we get the distinguished triangle
\begin{displaymath} \begin{CD}
\bb{K}_{Y} @> ad_{q'} >> Rq'_{*}\bb{K}_{E} @>>> i^{*}{\cal K} @> [1] >> \:,
\end{CD} \end{displaymath}
with $ad_{q'}: \bb{K}_{Y} \stackrel{\sim}{\to} R^{0}q'_{*}\bb{K}_{E}$ an isomorphism.
In particular 
\[{\cal H}^{k}(i^{*}{\cal K}) =
\begin{cases}
R^{k}q'_{*}(\bb{K}_{E}) = \bb{K}_{Y} &\text{for $k=2,\dots,2m-2$,}\\
0 &\text{otherwise.}
\end{cases} \]
But then $Hom({\cal K},\bb{K}_{X}[1])=0$ so that the triangle defining
${\cal K}$ {\em splits\/}:
\[Rq_{*}\bb{K}_{X'} \simeq \bb{K}_{X} \oplus {\cal K}\:.\]
Similarly
\[Hom\bigl(\bb{K}_{X}[n(X)] , \,^{m}{\cal H}^{0}({\cal K}[n(X')]) \bigr) = 0\]
so that 
\[(\,^{m}{\cal H}^{0}(Rq_{*}\bb{K}_{X'}[n(X')]),d) = 
(\bb{K}_{X}[n(X)],d) \oplus (\,^{m}{\cal H}^{0}({\cal K}[n(X')]),d') \]
for a suitable duality $d'$. 
Moreover this is already true for $\bb{K}=\bb{Z}$ so that
\[ [(\,^{m}{\cal H}^{0}(Rq_{*}\bb{K}_{X'}[n(X')]),d)] = 
[(\bb{K}_{X}[n(X)],d)] \:+\:  
\epsilon \cdot  i_{*}\bigl([(\bb{K}_{Y}[n(Y)],d)]\bigr) \]
for a suitable (locally constant) $\epsilon \in \{-1,0,1\}$. But then we can assume $\bb{K}=\bb{R}$
so that by Lemma \ref{lem:perverse} and \cite[cor.4.8]{CS1}:

\begin{equation} \label{eq:proof2}
q_{*}\bigl([(\bb{K}_{X'}[n(X')],d)]\bigr) = [(\bb{K}_{X}[n(X)],d)] \:+\: 
sign(E_{y})\cdot i_{*}\bigl([(\bb{K}_{Y}[n(Y)],d)]\bigr) \:.
\end{equation}

Here
\[ E_{y}:= q^{-1}N(y)/q^{-1}L(y) \simeq q^{-1}N(y)\cup_{q^{-1}L(y)}\; cone(q^{-1}L(y)) \]
for $y\in Y$ is obtained from the inverse image of a ``normal slice'' $N(y)$ to $Y$ at $y$
by collapsing the boundary to a point, or equivalently, by attaching the cone on the
boundary (compare \cite[p.522,thm.4.7]{CS1}). Note that the arguments of \cite{CS1}
work for our proper map $q$ without assuming $X'$ to be compact!

But in our case we can take 
\[y=pt=0 \in N(y)=(N_{Y}X)_{y} = \bb{C}^{m} \subset \bb{P}(\bb{C}^{m}\oplus \bb{C}) 
= \bb{P}^{m}(\bb{C}) \:,\]
with $L(y)\subset \bb{C}^{m}$ a sphere around $y=pt=0$. And to calculate $sign(E_{y})$,
we make the following trick: By ``Novikov additivity'' of  the {\em signature\/} 
(compare \cite[thm.(A)]{J} and \cite[prop.3.1]{Si}) we get
\[sign(E_{y})= sign(Bl_{\{pt\}}\bb{P}^{m}(\bb{C})) - sign(\bb{P}^{m}(\bb{C})) \:,\]
with $Bl_{\{pt\}}\bb{P}^{m}(\bb{C})$ the blow-up of $\bb{P}^{m}(\bb{C})$ along
the point $y=pt=0$. But by the {\em blow-up formula\/} for the signature, i.e.
equation (\ref{eq:blup}) with $y=1$, we also have
\[sign(E_{y})= sign(\bb{P}^{m-1}(\bb{C})) -1 \:.\]
So we finally get

\begin{equation} \label{eq:proof3}
\begin{split}
& q_{*}\bigl([(\bb{K}_{X'}[n(X')],d)]\bigr) = \\
 [(\bb{K}_{X}[n(X)],d)] \:&+\: 
(sign(\bb{P}^{m-1}(\bb{C})) -1)\cdot i_{*}\bigl([(\bb{K}_{Y}[n(Y)],d)]\bigr) \:,
\end{split} \end{equation}
which together with (\ref{eq:proof1}) implies the ``blow-up formula'' $(2.)$
of Corollary \ref{cor:bl}.
\end{proof}

Note that the normalization condition of Theorem \ref{thm:sd} implies
the commutativity of the diagram (\ref{eq:real2}). 

\begin{rem} \label{rem}
Note that our proof of Theorem \ref{thm:sd} does not apply
to the {\em Witt group\/} $W_{\pm}(X)$, instead of the cobordism group
$\Omega_{\pm}(X)$, since it depends on Lemma \ref{lem:perverse}.
\end{rem}

Let us now discuss
the behaviour of the transformation $sd$ under smooth pullback.
Let $f: X'\to X$ be a {\em smooth morphism\/} of complex algebraic varieties
(or compact complex spaces) of {\em constant relative dimension\/} $d(f)$.
Then one has a {\em duality isomorphism\/} (\cite[prop.3.3.2]{KS}, \cite[rem.4.2.3]{Sch1}):
\[f^{!}\simeq f^{*}[2d(f)] \]
and $f^{*}$ maps $D_{c}^{b}(X)$ to $D_{c}^{b}(X)$. By the canonical isomorphism
(\cite[prop.3.1.13]{KS}, \cite[cor.4.2.2]{Sch1})
\[f^{!}\circ D_{X} \simeq D_{X'}\circ f^{*}\]
we get 
\begin{equation} \label{eq:pullback}
(T^{d(f)}f^{*})\circ D_{X} \simeq (T^{-d(f)}f^{!})\circ D_{X} \simeq
(T^{-d(f)}D_{X'})\circ f^{*} \simeq D_{X'}\circ (T^{d(f)}f^{*}) 
\end{equation}
so that $T^{d(f)}f^{*}=f^{*}[d(f)]$ {\em commutes with duality\/}. And similarly
\[\begin{split}
f^{!}\circ (T^{2n(X)}D_{X}) &\simeq (T^{2n(X)}D_{X'})\circ f^{*} \\
& \simeq (T^{2n(X')}D_{X'})\circ (T^{2d(f)}f^{*}) \simeq (T^{2n(X')}D_{X'})\circ f^{!} \:,
\end{split} \]
with $n(X),n(X')$ the corresponding complex dimension of $X,X'$.
By our identification of twisted cobordism groups, this corresponds to
the following commutative diagram:
\begin{displaymath} \begin{CD}
\Omega(X') @> T^{n(X')} > \sim > \Omega(X',T^{2n(X')}D_{X'}) \\
@A f^{*}[d(f)] AA  @AA f^{!} A \\
\Omega(X) @> T^{n(X)} > \sim > \Omega(X,T^{2n(X)}D_{X}) \:.
\end{CD} \end{displaymath}

Note that $f^{*}[d(f)]$ is {\em parity\/} preserving (or changing) depending on $d(f)$ even (or odd).
If moreover $X$ (and therefore also $X'$) is smooth of pure dimension,
then one gets by definition
\[f^{*}[d(f)]\bigl([(\bb{K}_{X}[n(X)],d)]\bigr) = [(\bb{K}_{X'}[n(X')],d)] \:.\]
And this implies the

\begin{cor} \label{cor:sdpullback}
Let $f: X'\to X$ be a {\em smooth morphism\/} of complex algebraic varieties
(or compact complex spaces) of {\em constant relative dimension\/} $d(f)$.
Then the following diagram commutes:
\begin{displaymath} \begin{CD}
K_{0}(var/X) @> f^{*} >> K_{0}(var/X') \\
@V sd VV @VV sd V \\
\Omega(X) @>> f^{*}[d(f)] > \Omega(X') \:.
\end{CD} \end{displaymath}
If moreover $X$ and $X'$ are compact, then 
the following {\em Verdier Riemann-Roch type\/} diagram
\begin{displaymath} \begin{CD}
K_{0}(var/X) @> f^{*} >> K_{0}(var/X') \\
@V L_{*}\circ sd VV @VV L_{*}\circ sd V \\
H_{2*}(X,\bb{Q}) @>> L^{*}(T_{f})\cap f^{*} > H_{2*}(X',\bb{Q}) 
\end{CD} \end{displaymath} 
is also commutative. \hfill $\Box$
\end{cor}

\begin{rem} \label{rem:sdVRR}
We expect a similar {\em Verdier Riemann-Roch theorem\/} for the $L$-class transformation
$L_{*}$, but at the moment we have no proof or reference for it.
\end{rem}

Let us consider a {\em fixed\/} complex analytic (algebraic) Whitney stratification
$X_{\bullet}$ of a closed subset $X$ in a complex manifold $M$, with
$f: X':=X\cap M' \to X$ the inclusion of the intersection with a closed
complex submanifold $M'$ of $M$ which is {\em transversal\/} to $X_{\bullet}$.
Then $X'$ gets an induced complex Whitney stratification $X'_{\bullet}$.
Let $D^{b}_{c}(X_{\bullet})$ be the bounded derived category of sheaf complexes
which are constructible with respect to this stratification (and similarly
for $X'_{\bullet}$, compare \cite[chap.4]{Sch1}).
Then we get by \cite[cor.5.4.11]{KS} as before a commutative diagram:
\begin{displaymath} \begin{CD}
\Omega(X_{\bullet}') @> T^{n(X')} > \sim > \Omega(X_{\bullet}',T^{2n(X')}D_{X'}) \\
@A f^{*}[d(f)] AA  @AA f^{!} A \\
\Omega(X_{\bullet}) @> T^{n(X)} > \sim > \Omega(X_{\bullet},T^{2n(X)}D_{X}) \:,
\end{CD} \end{displaymath}
with $-d(f):=n(M)-n(M')$ the complex codimension of $M'$.\\

And again it is natural to ask for $X,X'$ compact a corresponding {\em Verdier Riemann-Roch theorem\/}
for the $L$-class transformation $L_{*}$, with 
\[T_{f}:=-N_{M'}M|X' \in K^{0}(X') \]
(compare with \cite{Sch2} for the corresponding result in the context of
the Chern-Schwartz-MacPherson transformation $c_{*}$).\\

At least in the case of a {\em trivial normal bundle\/} $N_{M'}M|X'$, this is true by
\cite[thm.5.1]{CS1}. In fact, this property for ``oriented stratified spaces'' together with
the description of the {\em degree\/} of $L_{0}$ as in Theorem \ref{thm:Lclass} 
characterizes the $L$-class transformation $L_{*}$ of Cappell-Shaneson {\em uniquely\/}
\cite[thm.5.1]{CS1}. But here it is important to go outside the realm of 
complex analytic  stratifications!\\

Let us finish this section with the {\em multiplicativity\/} of our transformations with respect to {\em exterior
products\/}. Consider two complex algebraic (or compact analytic) spaces $X$ and $X'$.
Then one has by \cite[cor.2.0.4]{Sch1} a natural isomorphism
\begin{equation} \label{eq:extdual}
(D_{X}{\cal F}) \boxtimes (D_{X'}{\cal F}')
\simeq D_{X\times X'}({\cal F} \boxtimes {\cal F}') 
\end{equation}
for ${\cal F}\in D^{b}_{c}(X)$ and ${\cal F}'\in D^{b}_{c}(X')$.
Assume now that $({\cal F},d)$ is also selfdual. Then 
\[{\cal F} \boxtimes (\cdot): D^{b}_{c}(X')\to D^{b}_{c}(X\times X')\]
{\em commutes with duality\/} by the isomorphism
\begin{displaymath} \begin{CD}
{\cal F}\boxtimes (D_{X'}{\cal F}') @> d\;\boxtimes id > \sim >
(D_{X}{\cal F})\boxtimes (D_{X'}{\cal F}')
\simeq D_{X\times X'}({\cal F} \boxtimes {\cal F}') \:.
\end{CD} \end{displaymath}
And similarly for a selfdual $({\cal F}',d')$. In this way we get a {\em bilinear 
pairing\/}:

\begin{equation} \label{eq:extcob}
\begin{split}
& \times\;:\Omega(X) \times  \Omega(X')\to \Omega(X \times X')\;:\\
& [({\cal F},d)]\times [({\cal F}',d')] \mapsto
[({\cal F} \boxtimes {\cal F}',d\;  \boxtimes d')] \:.
\end{split} \end{equation}

Here $({\cal F} \boxtimes {\cal F}',d\;  \boxtimes d')$ is {\em symmetric} if
$({\cal F},d)$ and $({\cal F}',d')$ are of the same parity.
Compare also with \cite{GN} for the corresponding result for {\em Witt groups\/}
of abstract triangulated categories with duality.\\

If moreover $X$ and $X'$ are smooth and pure dimensional,
then one gets by definition
\[[(\bb{K}_{X}[n(X)],d)]\times [(\bb{K}_{X'}[n(X')],d)] = 
[(\bb{K}_{X\times X'}[n(X\times X')],d)] \:.\]
And this implies the

\begin{cor} \label{cor:sdproduct}
Consider two complex algebraic (or compact analytic) spaces $X$ and $X'$.
Then the following diagram commutes:
\begin{displaymath} \begin{CD}
K_{0}(var/X) \times K_{0}(var/X') @> \times >> K_{0}(var/X\times X') \\
@V sd \times V sd  V @VV sd V \\
\Omega(X) \times \Omega(X') @>> \times > \Omega(X\times X') \:.
\end{CD} \end{displaymath}
If moreover $X$ and $X'$ are compact, then the following  diagram
\begin{displaymath} \begin{CD}
K_{0}(var/X) \times K_{0}(var/X') @> \times >> K_{0}(var/X\times X') \\
@V L_{*}\circ sd  \times V L_{*} \circ sd V @VV L_{*}\circ sd V \\
H_{2*}(X,\bb{Q}) \times H_{2*}(X',\bb{Q})  @>> \times > H_{2*}(X\times X',\bb{Q}) 
\end{CD} \end{displaymath} 
is also commutative. \hfill $\Box$
\end{cor}

\begin{rem} \label{rem:sdproduct}
We expect a similar ``multiplicativity'' for the {\em $L$-class transformation\/}
$L_{*}$, but at the moment we have no proof or reference for it.
\end{rem}

\section{Hodge theoretic definition of motivic Chern classes}
\label{sec:Hodge}

In this section we explain another construction of the {\em motivic Chern class transformation\/} 
$mC_{*}$ with the help of some fundamental results from the theory of
{\em algebraic mixed Hodge modules\/} due to M.Saito \cite{Sai1}-\cite{Sai6}.
This is the most functorial (but also the most difficult) approach to our motivic
classes. In fact this is the way we found them first!
Moreover, this functorial approach is needed, if one wants to extend
the ``generalized Verdier Riemann-Roch theorem'' and the theory of ``Milnor
classes for local complete intersections'' \cite{Sch2} from the context of Chern classes
to the context of our new motivic characteristic classes.\\

Since this theory of algebraic mixed Hodge modules is a very complicated
(and far away from geometry), we reduce our construction to a
few formal properties, together with a simple and instructive
explicite calculation for the normalization condition,
all of which are contained in the work of M.Saito.\\

Let us assume that our base field is $k=\bb{C}$.
To motivate the following constructions, let us first recall the
definition of the {\em Hodge characteristic\/} transformation
$Hc: K_{0}(var/pt)\to \bb{Z}[u,v] \:.$
By the now classical theory of Deligne \cite{De1,De2,Sr}, the cohomology groups
$V=H^{i}_{c}(X^{an},\bb{Q})$ of a complex algebraic variety
have a canonical  functorial {\em mixed Hodge structure\/},
which includes in particular the following data on the finite dimensional
rational vector space $V$:
\begin{itemize}
\item A finite increasing (weight) filtration $W$ of $V$ with
$W_{i}=\{0\}$ for $i<<0$ and $W_{i}=V$ for $i>>0$.
\item A finite decreasing (Hodge) filtration $F$ of $V\otimes \bb{C}$ with
$F^{p}=V$ for $p<<0$ and $F^{p}=\{0\}$ for $p>>0$.
\end{itemize}
These filtrations have to satisfy some additional properties, which imply
that the transformation of taking suitable graded vector spaces
$gr^{W}_{i}, gr_{F}^{p}$ and $gr_{F}^{p}gr_{i}^{W}$ for $i,p\in \bb{Z}$
induce corresponding transformations on the Grothendieck group
$K^{0}(MHS)$ of the abelian category of (rational) mixed Hodge structures,
i.e. morphism of mixed Hodge structures are ``strictly stable'' with respect
to the filtrations $F$ and $W$. Assume $Y$ is a closed algebraic subset of
$X$ with open complement $U:=X\backslash Y$. Then the maps in the long exact
cohomology sequence
\begin{equation} \label{eq:long}
\cdots \to H^{i}_{c}(U^{an},\bb{Q}) \to H^{i}_{c}(X^{an},\bb{Q}) \to
H^{i}_{c}(Y^{an},\bb{Q}) \to \cdots
\end{equation}
are morphisms of mixed Hodge structures so that the function
\begin{equation} \label{eq:addgr}
X\mapsto Hc(X):= \sum_{i,p,q\geq 0}\:(-1)^{i}(-1)^{p+q}\cdot 
dim_{C}\left(gr_{F}^{p}gr_{p+q}^{W}H^{i}_{c}(X^{an},\bb{C})\right) u^{p}v^{q}
\end{equation} 
satisfies the ``additivity property'' (add). In this way we get the 
{\em Hodge characteristic\/} (compare \cite{Sr}):
\[Hc: K_{0}(var/\{pt\})\to \bb{Z}[u,v],\: [X\to \{pt\}]\mapsto Hc(X).\]
Note that most references do not (!) include the sign-factor
$(-1)^{p+q}$ in their definition of the Hodge characteristic,
which is then called the {\em $E$-polynomial\/}:
\begin{equation} \label{eq:E}
E(X)(u,v) := Hc(X)(-u,-v) \:.
\end{equation}
Our sign convention fits better with the following normalization for
$X$ smooth and complete.
Specializing further, one gets also the (compare \cite{Lo})
\begin{itemize} \item
{\em Hodge filtration characteristic\/} corresponding to $(u,v)=(y,-1)$:
\[Hfc(X):=\sum_{i,p\geq 0}\:(-1)^{i} 
dim_{C}\left(gr_{F}^{p}H^{i}_{c}(X^{an},\bb{C})\right) (-y)^{p} \:.\] \item
{\em Weight filtration characteristic\/} corresponding to $(u,v)=(w,w)$: \[wc(X):=
\sum_{i,q\geq 0}\:(-1)^{i}\cdot 
dim_{C}\left(gr_{q}^{W}H^{i}_{c}(X^{an},\bb{Q})\right) (-w)^{q}\:.\] \item
{\em Euler characteristic\/} (with compact support)  corresponding to $(u,v)=(-1,-1)$:
\[e(X):= \sum_{i\geq 0}\:(-1)^{i}\cdot 
dim_{C}(H^{i}_{c}(X^{an},\bb{Q})) \:. \]
\end{itemize}
So these specializations fit into the following commutative diagram:
\begin{displaymath} \begin{CD}
\bb{Z}[u,v] @> u=y > v=-1 > \bb{Z}[y] \\
@V u=w V v=w V  @VV y=-1 V \\
\bb{Z}[w] @>> w=-1 > \bb{Z}\:.
\end{CD} \end{displaymath}
Finally, this classical Hodge theory \cite{De1,De2,Sr} implies for $X$ smooth complete
(of pure dimension $d$)
the ``purity result'' $gr_{p+q}^{W}H^{i}(X^{an},\bb{C})=0$ for $p+q\neq i$,
together with
\begin{align*}
h^{p,q}(X):= &\sum_{i\geq 0}\: (-1)^{i}(-1)^{p+q}\cdot
dim_{C}\left(gr_{F}^{p}gr_{p+q}^{W}H^{i}_{c}(X^{an},\bb{C})\right)\\
= & dim_{C}\left(gr_{F}^{p}H^{p+q}(X^{an},\bb{C})\right)
= dim_{C}H^{q}(X^{an},\Lambda^{p}T^{*}X^{an})\\
= &  dim_{C}H^{q}(X,\Lambda^{p}T^{*}X)\:.
\end{align*}
Here the last two equalities follow from GAGA, the degeneration of the
``Hodge to de Rham spectral sequence''
\[E^{p,q}_{1}= H^{q}(X^{an},\Lambda^{p}T^{*}X^{an}) \to
H^{p+q}(X^{an},\Lambda^{\bullet}T^{*}X^{an})\]
at $E_{1}$ and the ``holomorphic Poincar\'{e} lemma''
\[H^{p+q}(X^{an},\bb{C})\simeq H^{p+q}(X^{an},\Lambda^{\bullet}T^{*}X^{an})
\:.\] So the holomorphic de Rham complex 
$DR({\cal O}_{X^{an}}):=[\Lambda^{\bullet}T^{*}X^{an}]$ (with ${\cal
O}_{X^{an}}$ in degree zero) is a resolution of the constant sheaf $\bb{C}$ on
$X^{an}$, and the ``stupid decreasing filtration'' 
\begin{equation} \label{eq:stup}
F^{p}DR({\cal O}_{X^{an}}):=[0 \to \cdots \to 0 \to \Lambda^{p}T^{*}X^{an} \to
\cdots  \Lambda^{d}T^{*}X^{an}] 
\end{equation}
induces the Hodge filtration $F$ on $H^{*}(X^{an},\bb{C})$.\\

In particular $T_{y*}([X\to \{pt\}])= Hfc([X\to \{pt\}])$ for $X$ smooth and complete
by (gHRR). But these classes $[X\to \{pt\}]$ generate $K_{0}(var/\{pt\})$ so
that we get the following Hodge theoretic description for any $X$:
\begin{equation} \label{eq:TyH}
\chi_{y}(X)= T_{y*}([X\to \{pt\}]) = \sum_{i,p\geq 0}\:(-1)^{i} 
dim_{C}\left(gr_{F}^{p}H^{i}_{c}(X^{an},\bb{C})\right) (-y)^{p} \:.
\end{equation}

And exactly this description can be generalized to the context of
relative Grothendieck groups $K_{0}(var/X)$ using the machinery of
mixed Hodge modules of M.Saito. But before we explain this, let us
point out another remark. All our characteristics above are indeed
ring homomorphisms on $K_{0}(var/\{pt\})$, because this is the
case for $\chi_{y}=Hfc$ (for example by remark \ref{rem:local}).
Such ring homomorphism are called ``characteristics'' \cite{DL,Lo}
or sometimes also ``motivic measures'' \cite{LL},
and there are much more examples known. In this sense
our transformation $T_{y*}$ is certainly a {\em motivic characteristic 
class\/} since it is a homology class version of the
motivic characteristic $Hfc$, just like the Chern-Schwartz-MacPherson
(class) transformation $c_{*}$ is a homology class version
of the Euler-Poincar\'{e} characteristic $e$. 

\begin{rem} \label{rem:ring}
Our {\em motivic characteristic classes\/} are only group homomorphisms,
because the corresponding {\em homology theories\/} are only groups for
a singular space $X$. But they commute with exterior products,
so that they are ring homomorphisms for $X=\{pt\}$.
But this is in general no longer true for $X$ smooth,
even if one has then a corresponding ring structure on the
{\em (co)homology\/}. This is closely related to a corresponding
{\em Verdier Riemann-Roch formula\/} for the diagonal embedding 
$d: X\to X\times X$ (compare \cite{Sch2,Y5}).
\end{rem}

One can ask if also the other characteristic
$Hc$ or $wc$ can be ``lifted up'' to such a homology class
transformation. But here the answer will be {\em no\/} (cf. \cite{Jo})!

\begin{ex} \label{ex:counter}
Assume that there is a functorial transformation 
\[T_{u,v}: K_{0}(var/X)\to H_{*}(X)\otimes\bb{Q}[u,v]\]
commuting with proper
pushdown, and also with pullback $f^{*}$ for a finite smooth morphism
$f: X' \to X$ between smooth varieties, such that for $X=\{pt\}$ we get back the
{\em Hodge characteristic\/} $Hc$. Let $d$ be the degree of such a covering map
$f$. Then it follows (with $k: X\to \{pt\}$ a constant map):
\[Hc(X')= T_{u,v}([X'\to \{pt\}]) = T_{u,v}(k_{*}f_{*}f^{*}[id_{X}]) =
k_{*}f_{*}f^{*}T_{u,v}([id_{X}])\]
\[ = k_{*}(d\cdot T_{u,v}([id_{X}])) = d \cdot k_{*}T_{u,v}([id_{X}])
= d\cdot Hc(X) = Hc(d\cdot [id_{pt}]) \cdot Hc(X)  \:.\]
So (as usual), the transformation $Hc$ has then to be multiplicative
in such finite coverings. But this is not the case.
Let $X'\to X$ be such a finite covering of degree $d>1$ over an elliptic
curve $X$. Then $X'$ is also an elliptic curve so that 
\[  Hc(X) = Hc(X') = (1+u)(1+v) \neq 0\:.\]
Note that the same argument applies also to the {\em weight characteristic\/},
with $wc(X)=wc(X')=(1+w)^{2}\neq 0$. Of course, everything goes well
in the context of $e$ and $Hfc$, since both are zero for an elliptic 
curve!
\end{ex}

Let us now formulate those results about algebraic mixed Hodge modules,
which we need for our application to the motivic Chern class transformation
$mC_{*}$. All these results are contained in the deep and long work of M.Saito.
Since most readers will not be familiar with this theory,
we present them in an axiomatic way pointing out the similarities
to constructible functions $F(X)$ and motivic Grothendieck groups 
$K_{0}(var/X)$.\\

Let $k$ be a subfield of $\bb{C}$. Then we work in the category
of reduced seperated schemes of finite type over $spec(k)$,
which we also call ``spaces'' or ``varieties'', with $pt=spec(k)$.

\begin{enumerate}
\item[MHM1:] To such a space $X$ one can associate an abelian category
of {\em algebraic mixed Hodge modules\/} $MHM(X/k)$, together with a functorial
pullback $f^{*}$ and pushdown $f_{!}$ on the level of derived categories
$D^{b}MHM(X/k)$ for any (not necessarily proper) map \cite[sec.4]{Sai2}
(and compare also with \cite{Sai5,Sai6}).
These transformations are functors of triangulated categories.
\item[MHM2] Let $i:Y\to X$ be the inclusion of a closed subspace, 
with open complement $j: U:=X\backslash Y \to X$.
Then one has for $M\in D^{b}MHM(X/k)$ a distinguished triangle (\cite[(eq.(4.4.1),p.321]{Sai2})
\[j_{!}j^{*}M \to M \to i_{!}i^{*}M \stackrel{[1]}{\to} \:.\]
\item[MHM3:] For all $p\in \bb{Z}$ one has a functor of triangulated
categories
\[gr^{F}_{p}DR: D^{b}MHM(X/k) \to D^{b}_{coh}(X)\]
commuting with proper pushdown (compare with \cite[sec.2.3]{Sai1},
\cite[p.273]{Sai2}, \cite[eq.(1.3.4),p.9, prop.2.8]{Sai5} and also with 
\cite{Sai4,Sai5}). Here $D^{b}_{coh}(X)$ is
the bounded derived category of sheaves of ${\cal O}_{X}$-modules 
with coherent cohomology sheaves. Moreover, $gr^{F}_{p}DR(M)=0$
for almost all $p$ and $M\in D^{b}MHM(X/k)$ fixed (\cite[prop.2.2.10,
eq.(2.2.10.5)]{Sai1} and \cite[lem.1.14]{Sai3}).
\item[MHM4:] There is a distinguished element $\bb{Q}_{pt}^{H}
\in MHM(\{pt\}/k)$ such that 
\[gr^{F}_{-p}DR(\bb{Q}_{X}^{H}) \simeq \Lambda^{p}T^{*}X[-p]
\in D^{b}_{coh}(X)\]
for $X$ smooth and pure dimensional (\cite{Sai2}).
Here $\bb{Q}_{X}^{H}:=k^{*}\bb{Q}_{pt}^{H}$ for $k: X\to \{pt\}$
a constant map, with $\bb{Q}_{pt}^{H}$ viewed as a complex
concentrated in degree zero.
\end{enumerate}

Since the above transformations are functors of triangulated
categories, they induce functors on the level of {\em Grothendieck groups\/}
(of triangulated categories)
which we denote by the same name. By general nonsense one gets
for these {\em Grothendieck groups\/} isomorphisms
\[K_{0}(D^{b}MHM(X/k)) \simeq K_{0}(MHM(X/k)) \quad \text{and}
\quad K_{0}(D^{b}_{coh}(X)) \simeq G_{0}(X)\]
by associating to a complex its alternating sum of cohomology objects.

As explained later on, $K_{0}(MHM(X/k))$ plays the role of
``Hodge constructible functions'' with the class of $\bb{Q}_{X}^{H}$
as ``constant Hodge function'' on $X$.
By (MHM3) we get a group homomorphism commuting with proper pushdown:

\begin{equation} \label{eq:grH}
\begin{split}
 gr^{F}_{-*}DR: K_{0}(MHM(X/k)) \to G_{0}(X)\otimes \bb{Z}[y,y^{-1}] \:;\\
[M] \mapsto \sum_{p} \: [gr^{F}_{-p}DR(M)]\cdot (-y)^{p} \:.
\end{split} \end{equation}

And as for the map $e$ from motivic Grothendieck groups to constructible
functions, we get a group homomorphism commuting with pushdown
(compare also with \cite[sec.4]{Lo}):
\begin{equation} \label{eq:mH}
mH: K_{0}(var/X) \to K_{0}(MHM(X/k))\:, \:
[f: X'\to X] \mapsto [f_{!} \bb{Q}_{X'}^{H}]  \:.
\end{equation}

Indeed, the ``additivity relation'' (add) follows from (MHM2)
together with the functoriality of pushdown and pullback (MHM1):
For $i: Y\to X'$ the inclusion of a closed subspace, with open
complement $j:U\to X'$, the distinguished triangle
(with $k$ the constant map on $X'$)
\[j_{!}j^{*}k^{*}\bb{Q}_{pt}^{H} \to k^{*}\bb{Q}_{pt}^{H} 
\to i_{!}i^{*}k^{*}\bb{Q}_{pt}^{H} \stackrel{[1]}{\to} \]
induces under $f_{!}$ the distinguished triangle (with $f: X'\to X$ as before)
\[f_{!}j_{!}j^{*}k^{*}\bb{Q}_{pt}^{H}  \to f_{!}k^{*}\bb{Q}_{pt}^{H} 
\to f_{!}i_{!}i^{*}k^{*}\bb{Q}_{pt}^{H} \stackrel{[1]}{\to} \:.\]
It translates in the corresponding Grothendieck group into the relation
\[[f_{!}k^{*}\bb{Q}_{pt}^{H}] = [f_{!}j_{!}j^{*}k^{*}\bb{Q}_{pt}^{H}]
+ [f_{!}i_{!}i^{*}k^{*}\bb{Q}_{pt}^{H}]\:.\]
This finally is nothing else than the asked additivity property 
\[[f_{!} \bb{Q}_{X'}^{H}] = [(f\circ j)_{!}\bb{Q}_{U}^{H}] +
 [(f\circ i)_{!}\bb{Q}_{Y}^{H}] \in K_{0}(D^{b}MHM(X/k)) \:.\]
Moreover, $mH$ commutes with pushdown 
for a map $f: X'\to X$ again by functoriality:
\[mH(f_{!}[g:Y\to X']) = mH([f\circ g: Y\to X]) =
[(f\circ g)_{!}\bb{Q}_{Y}^{H}] \]
\[= [f_{!}g_{!}\bb{Q}_{Y}^{H}] = f_{!}[g_{!}\bb{Q}_{Y}^{H}] 
= f_{!}mH([g:Y\to X])\:.\]

By (MHM4) we get for $X$ smooth and pure dimensional:
\[ gr^{F}_{-*}DR\circ mH([id_{X}]) =  
\sum_{i=0}^{dim X} \; [\Lambda^{i} T^{*}X]\cdot y^{i}
\in G_{0}(X)\otimes \bb{Z}[y,y^{-1}] \:.\]

\begin{cor} \label{cor:mc=DRmH}
The {\em motivic Chern class transformation\/} $mC_{*}$ of theorem \ref{thm:mC}
is given as the composition 
\[mC_{*}= gr^{F}_{-*}DR\circ mH: K_{0}(var/X) \to G_{0}(X)\otimes \bb{Z}[y]
\subset G_{0}(X)\otimes \bb{Z}[y,y^{-1}]\:. \quad \Box\]
\end{cor}

\begin{rem} \label{rem:emb}
The definition of $MHM(X/k)$ and therefore also the transformations
$mH$ and  $gr^{F}_{-*}DR$ depend a priory on the {\em embedding\/} $k\subset \bb{C}$.
By the uniqueness statement of Theorem \ref{thm:mC} and \ref{thm:Ty}
this is not the case for the transformations $mC_{*}$ and $T_{y*}$,
i.e. they are {\em independent\/} of the choice of the embedding $k\subset \bb{C}$,
if their definition is based on Corollary \ref{cor:mc=DRmH}.
\end{rem}

Let us now explain a little bit of the definition of the abelian category
$MHM(X/k)$ of algebraic mixed Hodge modules on $X/k$.
Its objects are special tuples $((M,F),K,W)$, which for $X$ smooth are given by
\begin{itemize}
\item $(M,F)$ an {\em algebraic holonomic filtered D-module\/} $M$ on $X$ with an
exhaustive,  bounded from below and increasing (Hodge) filtration $F$ by
algebraic ${\cal O}_{X}$-modules such that $gr_{*}^{F}M$ is a {\em coherent\/}
$gr_{*}^{F}{\cal D}_{X}$-module. In particular, the filtration $F$ is
 finite, which will imply the last claim of (MHM3).  Here the filtration $F$ on
the sheaf of algebraic differential operators ${\cal D}_{X}$ on $X$ is the
order filtration, and one can work either with left or right D-modules.
For singular $X$ one works with suitable local embeddings into manifolds
and corresponding filtered D-modules with support on $X$
(compare \cite{Sai1,Sai2,Sai4}).
\item $K\in D^{b}_{c}(X(\bb{C})^{an},\bb{Q})$ is an {\em algebraically constructible
complex of sheaves\/} of $\bb{Q}$-vector spaces (with finite dimensional stalks,
compare for example with \cite{Sch1}))
on the associated analytic space $X(\bb{C})^{an}$ corresponding to the induced
algebraic variety $X(\bb{C}):=X\otimes_{k}\bb{C}$ over $\bb{C}$,
which is {\em perverse\/} with respect to the middle perversity $t$-structure.
$F$ is called the underlying rational sheaf complex.
\item In addition one fixes a quasi-isomorphism $\alpha$ between 
$K(\bb{C}):=K \otimes_{Q}\bb{C}$ and the {\em holomorphic de Rham complex\/}
$DR(M(\bb{C})^{an})$ associated to the induced ${\cal D}_{X(C)}$-module
$M(\bb{C}):=M \otimes_{k}\bb{C}$.
\item $W$ is finally a finite increasing {\em (weight) filtration\/} of
$(M,F)$ and $K$, compatible in the obvious sense with the quasi-isomorphism
$\alpha$ above.
\end{itemize}
These data have to satisfy a long list of properties which we do not recall
here (since it is not important for us). In particular, one gets the
{\em equivalence\/} \cite[eq.(4.2.12),p.319]{Sai2})
\begin{equation} \label{eq:MHM=MHS}
\begin{split}
MHM(\{pt\}/\bb{C}) &\simeq \{\text{(graded) polarizable mixed $\bb{Q}$-Hodge
structures} \} \\
\:F^{-p} &\leftrightarrow F_{p} \end{split}
\end{equation}
between the category of {\em algebraic mixed Hodge modules on\/} $pt=spec(\bb{C})$,
and the category of {\em (graded) polarizable mixed $\bb{Q}$-Hodge
structures\/}. Of course, one has to switch the increasing $D$-module
filtration $F^{p}$ to the decreasing Hodge filtration by
$F^{-p}\leftrightarrow F_{p}$  so that $gr^{F}_{-p}\simeq gr_{F}^{p}$.
For elements in $MHM(\{pt\}/k)$, the corresponding Hodge filtration is already
defined over $k$ (compare \cite[sec.1.3]{Sai6}). 

The {\em distinguished element\/} $\bb{Q}^{H}_{pt}\in MHM(\{pt\}/k)$ of (MHM4) is given by
\begin{equation} \label{eq:disting}
\bb{Q}^{H}_{pt}:=((k,F),\bb{Q},W) \quad \text{with $gr^{F}_{i}=0=gr^{W}_{i}$
for $i\neq 0$}
\end{equation}
and $\alpha: k \otimes \bb{C}\simeq \bb{Q}\otimes \bb{C}$ the obvious
isomorphism (compare \cite[sec.1.3]{Sai6}).
The functorial pullback and pushdown of (MHM1) corresponds under the {\em forget
functor\/} (\cite[thm.0.1,p.222]{Sai2})
\begin{equation} \label{eq:rat}
rat: D^{b}MHM(X/k)\to D^{b}_{c}(X(\bb{C})^{an},\bb{Q}) \:,
((M,F),K,W)\mapsto K 
\end{equation}
to the classical corresponding (derived) functors $f^{*}$ and $f_{!}$
on the level of algebraically constructible sheaf complexes,
with $rat(\bb{Q}^{H}_{X})\simeq \bb{Q}_{X(C)^{an}}$.
So by (\ref{eq:MHM=MHS}) one should think of an algebraic mixed Hodge module
as a kind of ``(perverse) constructible Hodge sheaf''!
But one has to be very careful with this analogy. $\bb{Q}^{H}_{X}$
is in general a highly complicated complex in $D^{b}MHM(X/k)$,
which is impossible to calculate explicitely. But if $X$ is {\em smooth and
pure $d$-dimensional\/}, then $\bb{Q}_{X(C)^{an}}[d]$ is a {\em perverse sheaf\/}
and $\bb{Q}^{H}_{X}[d]\in MHM(X/k)$ a single {\em mixed Hodge module\/}
(in degree $0$), which is explicitely given by (\cite[eq.(4.4.2),p.322]{Sai2}):
\begin{equation} \label{MHMsmooth}
\bb{Q}^{H}_{X}[d] \simeq (({\cal O}_{X},F),\bb{Q}_{X(C)^{an}}[d],W) \:,
\end{equation}
with $F$ and $W$ the trivial filtration $gr_{i}^{F}=0=gr_{d+i}^{W}$
for $i\neq 0$. Here we use for the underlying D-module the description as the
left D-module ${\cal O}_{X}$, which maybe is more natural at this point.\\

The distingished triangle (MHM2) is a ``lift'' of the corresponding
distinguished triangle for constructible sheaves.
Similarly, by taking a constant map $f: X\to pt$ we get by (MHM1) and
(\ref{eq:MHM=MHS}) a functorial (rational) {\em mixed Hodge structure\/} on
\[rat \left(R^{i}f_{!}\bb{Q}^{H}_{X} \right) \simeq
H^{i}_{c}(X(\bb{C})^{an},\bb{Q}) \:,\]
whose Hodge numbers are easily seen to be the same as those coming
from the mixed Hodge structure of Deligne \cite{De1,De2,Sr} (both have the same additivity
property so that one only has to compare the case $X$ smooth, which
follows from the constructions). In fact, even the Hodge structures 
are the same by a deep theorem of M.Saito \cite[cor.4.3]{Sai5}.\\

Let us finally explain (MHM3) and (MHM4) for the case $X$ smooth and
pure $d$-dimensional. The {\em de Rham functor\/} $DR$ factorizes as
(compare \cite{Sai1,Sai3,Sai5})
\begin{equation} \label{eq:facDR}
DR: D^{b}MHM(X/k) \to D^{b}F_{coh}(X,Diff) \to  
D^{b}_{c}(X(\bb{C})^{an},\bb{C}) \:,
\end{equation}
with $D^{b}F_{coh}(X,Diff)$ the ``bounded derived category of
filtered differential complexes'' on $X$ with coherent cohomology sheaves.
Here the objects are bounded complexes $(L^{\bullet},F)$ of ${\cal
O}_{X}$-sheaves with an increasing (bounded from below) filtration $F$ by such
sheaves, whose morphisms are ``differential operators'' in a suitable sense. In
particular \[gr^{F}_{p}(L^{\bullet}) \in D^{b}(X,{\cal O}_{X})\]
becomes an ${\cal O}_{X}$-linear complex with coherent cohomology.
Moreover, the morphisms of mixed Hodge modules are ``strict'' with respect to the
Hodge filtration $F$ (and the weight filtration $W$) so that 
$gr^{F}_{p}DR$ induces the corresponding transformation of (MHM3).
Finally, $DR(\bb{Q}^{H}_{X})$ is given by the usual {\em de Rham complex\/}
$\Lambda^{\bullet}T^{*}X$ with the induced increasing filtration  
\begin{equation} \label{eq:DR2}
F_{p}DR(\bb{Q}^{H}_{X}) := [ F_{p}{\cal O}_{X} \to 
F_{p+1}{\cal O}_{X}\otimes \Lambda^{1} T^{*}X \to \: \cdots \to  
F_{p+d}{\cal O}_{X}\otimes \Lambda^{d} T^{*}X ]\:,
\end{equation}
with $F_{p}{\cal O}_{X}$ in degree zero and $F$ the trivial filtration with
$gr_{i}^{F}=0$ for $i\neq 0$. Let us switch to the corresponding decreasing
filtration (with $gr^{F}_{-p}\simeq gr^{p}_{F}$):
\[F^{p}DR(\bb{Q}^{H}_{X}) := [F^{p}{\cal O}_{X} \to 
F^{p-1}{\cal O}_{X}\otimes \Lambda^{1} T^{*}X \to \: \cdots \to  
F^{p-d}{\cal O}_{X}\otimes \Lambda^{d} T^{*}X] \:.\]
Then the {\em de Rham complex\/} $DR({\cal O}_{X})$ with the stupid filtration $\sigma^{p}$ as in (\ref{eq:stup}): 
\[F^{p}DR({\cal O}_{X}):= [0 \to \cdots \to 0 \to \Lambda^{p} T^{*}X \to
\cdots  \Lambda^{d} T^{*}X] \:, \]
becomes a filtered subcomplex.
And one trivially checks that the inclusion induces on the associated graded
complexes the isomorphism 
\[gr^{p}_{\sigma} DR({\cal O}_{X}) \simeq \Lambda^{p}T^{*}X[-p] 
 \simeq gr^{p}_{F}DR (\bb{Q}^{H}_{X}) 
\simeq gr_{-p}^{F}DR (\bb{Q}^{H}_{X}) \:.\]
In this way one finally also gets (MHM4). 

\begin{rem} \label{rem:totaro}
The use of the transformation $gr^{F}_{p}DR$ of (MHM3) in the context
of {\em characteristic classes of singular spaces\/} is not new.
It was already used by Totaro \cite{To} in his study of the relation between
{\em Chern numbers\/} for singular complex varieties and {\em elliptic homology\/}!

But he was interested in characteristic numbers and classes
invariant under {\em small resolution\/}, and not in {\em functoriality\/} as in
our paper. So  he works with the counterpart ${\cal IC}_{X}^{H}
\in MHM(X/\bb{C})$ of the {\em intersection cohomology complex\/}
instead of the constant ``Hodge sheaf'' 
$\bb{Q}_{X}^{H} \in D^{b}MHM(X/\bb{C})$ as used in this paper.
He then also applied the {\em singular Riemann-Roch transformation\/}
$td_{*}$ of Baum-Fulton-MacPherson to associate to a singular
complex algebraic variety $X$ of dimension $n$ some natural {\em homology classes\/}
$\chi^{n-k}_{p}(X)\in H^{BM}_{2k}(X,\bb{Q})$ for $p\in \bb{Z}$.
In our notation, the corresponding {\em total homology class\/}
$\chi^{n-*}_{p}(X)\in H^{BM}_{2*}(X,\bb{Q})$ 
is given by evaluating
\[td_{(1+y)}\circ gr^{F}_{-p}DR({\cal IC}_{X}^{H}) \in 
H^{BM}_{2*}(X,\bb{Q})[y,(1+y)^{-1}] \]
at $y=0$. Here it is important to work with the transformation 
\[td_{(1+y)}\circ gr^{F}_{-p}DR: K_{0}(MHM(X/\bb{C}))\to 
H^{BM}_{2*}(X,\bb{Q})[y,y^{-1},(1+y)^{-1}] \:.\]
This allows one to use more general
coefficients like ${\cal IC}_{X}^{H} \in MHM(X/\bb{C})$,
which are a priori not in the image of $mH$.
\end{rem}

\begin{rem}
The announcement \cite{CS2} and \cite[sec.4]{Sh} suggests for a 
pure $d$-dimensional compact complex algebraic
variety $X$ the following relation between the {\em Hodge theoretical classes\/}
and the {\em topological $L$-class\/} 
$L_{*}(X)=L_{*}({\cal IC}_{X})$ of Goresky-MacPherson:
\begin{equation} \label{eq:IC}
(\:td_{(1+y)}\circ gr^{F}_{-*}DR({\cal IC}_{X}^{H}[-d])\:)|_{y=1} = L_{*}(X) \:.
\end{equation}
At least the equality of their {\em degrees\/} follows from the work of Saito,
i.e. the description of the {\em signature\/} of the global Intersection (co)homology
in terms of {\em Hodge numbers\/} (as in the {\em Hodge index theorem\/} for smooth K\"ahler manifolds).

Note that one has a natural morphism $\bb{Q}_{X}^{H} \to
{\cal IC}_{X}^{H}[-d]$ in $D^{b}MHM(X/\bb{C})$, which is an
isomorphism for $X$ a {\em rational homology manifold\/}.
So (\ref{eq:IC}) would imply our {\em conjecture\/} 
that for a rational homology manifold $T_{1*}(X)=L_{*}(X)$,
and more generally it would explain the difference between
$T_{1*}(X)$ and $L_{*}(X)$. Similarly, (\ref{eq:IC}) implies that
$L_{*}(X)$ is in the image of the cycle map from the Chow group
$A_{*}(X)$ to homology.
\end{rem}

As explained before, it is in general impossible to calculate
$gr_{-p}^{F}DR (\bb{Q}^{H}_{X})$ explicitely for a singular
space $X$. But by comparing our different definitions of
$mC_{*}(X)=mC_{*}([id_{X}])$ in terms of the {\em Du Bois complex\/}
and in terms of {\em mixed Hodge modules\/}, we get a least (for all
$p\in \bb{Z}$):

\begin{equation} \label{eq:DBHodge}
[gr_{-p}^{F}DR (\bb{Q}^{H}_{X})] = 
[gr_{F}^{p}(\underline{\Omega}^{*}_{X})] \in G_{0}(X) \:.
\end{equation}

In the work of M.Saito \cite{Sai5} the reader will find 
a deeper identification of the underlying filtered complexes,
which is of course much stronger than the equality above on
the level of elements in the Grothendieck group.
Here we finally state only the following results from \cite{Sai5}
for a complex algebraic variety $X$ of dimension $n$:
\begin{enumerate}
\item $gr_{-p}^{F}DR (\bb{Q}^{H}_{X})\simeq 0 \in D^{b}_{coh}(X)$
for $p<0$ and $p>n$. 
\item Let $X'$ be a resolution of singularities of the union
of the $n$-dimensional irreducible components of $X$, with
$\pi: X'\to X$ the induced proper map. Then
\[ gr_{-n}^{F}DR (\bb{Q}^{H}_{X})\simeq \pi_{*} 
\Lambda^{n}T^{*}X'[-n] \:.\]
Note that $R^{i}\pi_{*} \Lambda^{n}T^{*}X'=0$ for $i>0$ by the
{\em Grauert-Riemenschneider vanishing theorem\/}. 
\item $h^{i}(gr_{0}^{F}DR (\bb{Q}^{H}_{X}))\simeq 0$ for $i<0$,
and 
\[h^{0}(gr_{0}^{F}DR (\bb{Q}^{H}_{X}))\simeq {\cal O}_{X}^{wn} \:,\]
with ${\cal O}_{X}^{wn}$ the coherent structure sheaf 
of the {\em weak normalization\/} $X_{max}$ of $X$ (whose underlying space
is identified with $X$). One gets in particular natural morphisms
\[{\cal O}_{X} \to {\cal O}_{X}^{wn} \to
gr_{0}^{F}DR (\bb{Q}^{H}_{X})\]
in $D^{b}_{coh}(X)$. And in this language $X$ has at most ``Du Bois singularities'' if 
the composed map ${\cal O}_{X} \to gr_{0}^{F}DR (\bb{Q}^{H}_{X})$ is a quasi-isomorphism. 
\end{enumerate}

$ $\\
Jean-Paul Brasselet\\
Institute de Math\'{e}matiques de Luminy\\
UPR 9016-CNRS\\
Campus de Luminy - Case 907\\
13288 Marseille Cedex 9, France\\
E-mail: jpb@iml.univ-mrs.fr\\
$ $\\
J\"{o}rg Sch\"{u}rmann\\
Westf. Wilhelms-Universit\"{a}t\\
SFB 478
"Geometrische Strukturen in der Mathematik" \\
Hittorfstr.27\\
48149 M\"{u}nster, Germany\\
E-mail: jschuerm@math.uni-muenster.de\\
$ $\\
Shoji Yokura\\
Department of Mathematics and Computer Science\\
Faculty of Science, University of Kagoshima\\
21-35 Korimoto 1-Chome\\
Kagoshima 890-0065, Japan\\
E-mail: yokura@sci.kagoshima-u.ac.jp

\end{document}